\documentclass[12pt, a4paper]{amsart}
\usepackage{amssymb,epsfig,color}

\newcommand{\be}{{\mathbf e}}

\newcommand{\bw}{{\mathbf w}}

\newcommand{\GK}{{G\!K}}
\newcommand{\ADO}{{A\!DO}}

\newcommand{\sresq}{\widetilde{\mathcal{U}}_q(sl_2)}
\newcommand{\Vol}{\operatorname{Vol}}
\newcommand{\CS}{\operatorname{CS}}
\newcommand{\Gk}{\operatorname{GK}}
\newcommand{\End}{\operatorname{End}}
\definecolor{gray}{rgb}{0.5,0.5,0.5}
\begin{document}
\title[Generalized Kashaev invariants for knots in three manifolds]{Generalized Kashaev invariants for knots in three manifolds}
\author[J. Murakami]{Jun Murakami}
\address{Department of Mathematics\\ Faculty of Science and Engineering\\
Waseda~University\\
3-4-1 Ohkubo, Shinjuku-ku \\
Tokyo~169-8555, Japan}
\email{murakami@waseda.jp}
\date{\today \\ \quad\  This work was supported in part by JSPS KAKENHI Grant Numbers 22540236, 25287014.}
\begin{abstract}Kashaev's invariants for a knot in a three sphere are generalized to invariants of a knot in a three manifold.  
A relation between the newly constructed invariants and the hyperbolic volume of the knot complement is observed  for some knots in lens spaces.
\end{abstract}
\maketitle
\baselineskip18pt
\section*{Introduction}
The Jones polynomial of knots and links is discovered in  \cite{Jo}, which is defined by a simple skein relation, and relates to the quantum enveloping algebra ${\mathcal U}_q(sl_2)$ through the quantum $R$-matrix.  
After the Jones polynomial, 
a large number of quantum invariants are constructed from various $R$-matrices associated with quantum enveloping algebras, Hopf algebras, and operator algebras.  
The Jones polynomial is also extended to invariants of three manifolds and links in three manifolds by \cite{Wi} and \cite{ReTu}.  
\par
On the other hand, from a study of quantum dilogarithm, 
R. Kashaev introduced an invariant of links in three manifolds in \cite{Ka4}.  
He also gave an $R$-matrix formulation of his invariants for knots in $S^3$, 
and  found in \cite{Ka} a relation between his invariants and the hyperbolic volumes of knot complements.        
Let $\left<K\right>_N$ be the Kashaev's invariant of a knot $K$ for a positive inveter $N$, then the relation he found is the following.    
\medskip
\par\noindent
{\bf Conjecture 1. (Kashaev's conjecture)}
\begin{it}
For a hyperbolic knot $K$ in $S^3$, 
$$
2\, \pi \, \lim_{N\to \infty} \dfrac{\log|\left<K\right>_N|}{N}
=
\Vol(K),   
$$
where $\Vol(K)$ is the hyperbolic volume of the knot complement $S^3 \setminus K$.  
\end{it}
\smallskip
\par
Kahsaev's invariant turned out to be a specialization of the colored Jones invariant in \cite{MuMu}, and the above conjecture is refined in \cite{MMOTY}  as follows.  
\medskip\par\noindent
{\bf Conjecture 2. (Complexification of Kashaev's conjecture)}
\begin{it}
For a hyperbolic knot $K$ in $S^3$, 
$$
\left<K\right>_N
\sim
\exp\dfrac{N}{2\, \pi}\left(\Vol(K) + \sqrt{-1} \, \operatorname{CS}(K)\right)
\qquad(N \to \infty)
$$
{\it where $\operatorname{CS}(K)$ is the Chern-Simons invariant \cite{CS}, \cite{Mey} of the knot complement  $S^3 \setminus K$.  }
\end{it}
\medskip\par
The above conjectures are not proved rigorously yet, but a method to obtain the hyperbolic volume and the Chern-Simons invariant from Kashaev's invariants are established in \cite{CMY} and \cite{Y}.  
\par
The aim of this paper is to construct certain quantum invariants for knots in three manifolds which have a relation to the hyperbolic volume as the above conjectures.  
We already have many quantum invariants for knots in three manifolds.  
Besides the invariants stated above, such invariants are constructed in \cite{Geer} and \cite{CGP} from finite-dimensional representations of the quantum group ${\mathcal U}_q(sl_2)$ at root of unity, and in \cite{KaLuVa} from the infinite dimensional representations of ${\mathcal U}_q(sl_2)$.  
However, it is not known about the actual relation between the above invariants and the hyperbolic volume of the complement of the knots.  
\par
Here we construct a family of invariants of a knot $\widetilde K$ in a three manifold $M$ by combining the Hennings invariant \cite{He} of three manifolds and the logarithmic invariant \cite{MN} of knots in $S^3$.   
This family contains a generalized Kashaev invariant $\Gk_N(\widetilde K)$, which coincides with Kashaev's invariant $\langle\widetilde K\rangle_N$ if $M = S^3$.  
Moreover, we introduce $\Gk_N^{SO(3)}(\widetilde K)$, which is the $SO(3)$ version of $\Gk_N(\widetilde K)$,  and  propose the following conjecture.  
\medskip\par\noindent
{\bf Conjecture 3. (Volume conjecture for the generalized Kashaev invariant)}
Let $\widetilde K$ be a knot in a three manifold $M$ such that the complement $M \setminus \widetilde K$ has the hyperbolic structure.  
Then
$$
\Gk_N^{SO(3)}(\widetilde K)
\sim
\exp\dfrac{N}{2\, \pi}\left(\Vol(\widetilde K) + \sqrt{-1} \, \operatorname{CS}(\widetilde K)\right)
\qquad(N \to \infty)  
$$
{\it where $\Vol(\widetilde K)$ and $\operatorname{CS}(\widetilde K)$ is the hyperbolic volume and the Chern-Simons invariant of the complement $M \setminus \widetilde K$.  }
\medskip\par\noindent
We give some examples for this conjecture at the end of this paper.  
\medskip\par\noindent
{\bf Remark 1.}
The invariants $GK_N(\tilde K)$ and $GK_N^{SO(3)}(\tilde K)$ are generalizations of Kashaev's invariant for knots in $S^3$.  
So they may have some relation to Kashaev's original invariant for knots in three manifolds defined in \cite{Ka4}.  
But any relation is not observed yet.  
\medskip\par
As we stated before, we construct invariants of knots in a three manifolds by combining the Hennings invariant and the logarithmic invariant.
The both of these invariants are related to the universal invariant 
introduced by Lawrence \cite{La} and Ohtsuki \cite{Oh1}, 
whose value is in a certain quotient of the small quantum group $\overline{\mathcal U}_q(sl_2)$,
which is  a finite dimensional Hopf algebra and is a quotient of the quantized enveloping algebra ${\mathcal U}_q(sl_2)$ where
 $q = e^{\pi i/N}$.
The generators and relations of  $\overline{\mathcal U}_q(sl_2)$ are given as follows.  
\begin{multline*}
\overline{\mathcal U}_q(sl_2)
=
\left<
K,\ K^{-1},\ E,\ F \mid
K \, E \, K^{-1} = q^2 \, E,\ 
K \, F \, K^{-1} =q^{-2} \, F,\right.\\ 
[E, F] = \dfrac{K-K^{-1}}{q-q^{-1}}, \  
\left.
\vphantom{\dfrac{K-K^{-1}}{q-q^{-1}}}
E^N = F^N = 0, \ K^{2N}=1
\right>.  
\end{multline*}
The Hopf algebra structure of $\overline{\mathcal U}_q(sl_2)$ is
given by 
the coproduct $\Delta : \overline{\mathcal U}_q(sl_2) \to \overline{\mathcal U}_q(sl_2) \otimes \overline{\mathcal U}_q(sl_2)$, the counit $\varepsilon : \overline{\mathcal U}_q(sl_2) \to \mathbf C$ and the antipode $S : \overline{\mathcal U}_q(sl_2) \to \overline{\mathcal U}_q(sl_2)$ satisfying 
$$
\begin{aligned}
\Delta(K) = K \otimes K, \ \Delta(E) = 1 \otimes E &+ E \otimes K,
\ \Delta(F) = K^{-1}\otimes F + F \otimes 1
\\
\epsilon(K)=1, \quad \epsilon(E)&=\epsilon(F)=0,
\\
S(K) = K^{-1}, \quad
S(E)&= -E\, K^{-1}, \quad
S(F) = -K\, F.
\end{aligned}
$$
The dimension of $\overline{\mathcal U}_q(sl_2)$ is $2 N^3$ and 
$$
\{E^a \, F^b \, K^c \mid 0 \leq a, \ b \leq N-1,\ 0 \leq c \leq 2\,N-1\}
$$ 
is a basis of it.  
\par
The universal invariant takes its value in the quotient $\overline{\mathcal U}_q(sl_2)/I$ where $I$ is the two sided ideal generated by commutators of $\overline{\mathcal U}_q(sl_2)$, i.e. 
$$
I = \left[\overline{\mathcal U}_q(sl_2), \overline{\mathcal U}_q(sl_2)\right] 
=
\left(x \, y - y \, x; \ x, y \in \overline{\mathcal U}_q(sl_2)\right).
$$  
\par
The Hennings invariant $H(M)$ for an oriented closed three manifold $M$ is constructed by using the {\it right integral} $\mu$, 
which is a linear functional on $\overline{\mathcal U}_q(sl_2)$ satisfying   
\begin{equation}
(\mu\otimes id)\Delta(x) = \mu(x) \, 1,
\label{eq:right integral}
\end{equation}
where $1$ is the unit of $\overline{\mathcal U}_q(sl_2)$.  
Such functional exists uniquely up to a scalar multiple since $\overline{\mathcal U}_q(sl_2)$  is a finite dimensional Hopf algebra.  
The relation \eqref{eq:right integral} corresponds to the second Kirby move and it allows  us to construct a three manifold invariant by using the right integral, which is the Hennings invariant.      
Let $\tau_N(M)$ be the Witten-Reshetihin-Turaev (WRT) invariant \cite{Wi}, \cite{ReTu} of $M$.  
Then it is shown in \cite{ChKuSr} and \cite{ChYuZh} that
the Hennings invariant can be expressed in terms of the WRT invariant for almost all cases.  
\par
Nagatomo and the author constructed in \cite{MN} the logarithmic invariant of a knot $K$ in $S^3$.  
Let ${\mathcal Z}$ be the center of  $\overline{\mathcal U}_q(sl_2)$.  
We study the center $c(T) \in {\mathcal Z}$ which corresponds to a tangle $T$ obtained from $K$, and we define knot invariants as the coefficients
of $c(T)$ with respect to certain basis of ${\mathcal Z}$.  
A topological quantum field theory (TQFT) based on the center ${\mathcal Z}$ is constructed by Kerler \cite{Ke}, and is refined by Feigin-Gainutdinov-Semikhatov-Tipunin \cite{FeGaSeTi} by using the logarithmic conformal field theory.    
The logarithmic knot invariant corresponds to this TQFT.   
\par
We also showed in \cite{MN} that the logarithmic invariant is expressed as a limit of the colored Alexander invariant, which is defined by Akutsu-Deguchi-Ohtsuki \cite{ADO} and is restudied by the author in \cite{JM}.    
It is an invariant of links with colored components, where the colors are complex numbers except integers.  
The logarithmic invariant is obtained as a limit of a sum of two colored Alexander invariants by taking its colors to certain integers.  
A relation like Conjectures 1, 2 are observed in \cite{CM} between the colored Alexander invariant and the hyperbolic volume of cone manifolds.  
\par
Let $M$ be a three manifold given by the surgery along a framed link $L$ in $S^3$, 
$\widetilde K$ be a knot in $M$, and  
$\widehat K$ be the pre-image of $\widetilde K$ in $S^3$.  
Then, to construct an invariant of $\widetilde K$, we apply the logarithmic invariant to $\widehat K$,  and apply the Henning invariant to $L$.  
\par
In Section 1, we recall the construction of the Hennings invariant and extend it to invariants of knots in three manifolds.  
In Section 2, we review irreducible and indecomposable representations of $\overline{\mathcal U}_q(sl_2)$.
By using these representations, we describe  centers and symmetric linear functions of $\overline{\mathcal U}_q(sl_2)$. 
In Section 3, we generalize the logarithmic invariants of knots in $S^3$ to invariants of knots in three manifolds. 
This family of invariants include the generalized Kashaev invariant $GK_N$.  
In Section 4, we investigate the generalized Kashaev invariant by using its relation to the colored Alexander invariant.  
In Section 5, we observe the relation between the generalized Kashaev invariants  of certain knots in lens spaces and the hyperbolic volumes of their  complements by numerical computation.  
\section{Colored Hennings invariants}
In this section, we generalize the colored invariants constructed by Hennings \cite{He}  for knots and links in $S^3$ to invariants for knots in three manifold, which we call the {\it colored Hennings invariant}.  
To do this, we first recall the construction of the universal  $\overline{\mathcal U}_q(sl_2)$ invariant for a link in $S^3$ introduced in \cite{La} and \cite{Oh1}.  
Then we apply Hennings' idea in \cite{He} to obtain invariants equipped with a color at each component of the link, where the color is given by a pair of a symmetric linear function and a center of $\overline{\mathcal U}_q(sl_2)$.  
There is a special symmetric linear function $\phi$ corresponding to the right integral $\mu$,  which assures the compatibility of the $\phi$ colored component with the second Kirby move.       
By using $\phi$, we construct invariants of links in arbitrary oriented three manifolds.  
If the  knot is emply, then this invariant coincides with the Hennings invariant of the three manifold introduced in \cite{He}, \cite{KaRa}, and \cite{Oh2} associated  with  $\overline{\mathcal U}_q(sl_2)$.  
\subsection{Notations}
Throughout this paper, let  $q = e^{\pi i/N}$.  
We use the following notations.  
$$
\{k\} = q^k -q^{-k}, 
\quad
\{k\}_+ = q^k + q^{-k}, \quad
[k] = \dfrac{\{k\}}{\{1\}}, 
\quad
[k]! = [k][k-1]\cdots [1], 
$$
$$
\{k\}! = \{k\}\{k-1\}\cdots \{1\}\ \  \text{for a positive integer $k$}, 
\quad
\{0\}! = [0]! = 1.  
$$
\subsection{The right-integral}
The {\it right integral} of a Hopf algebra is a non-trivial linear functional $\mu$ on the Hopf algebra which satisfies
$$
(\mu\otimes id)\Delta(x) = \mu(x)\, 1.  
$$
Any finite dimensional Hopf algebra has a right integral which is unique up to nonzero scalar multiplication.  
For detail, see \cite{Ra}.  
For $\overline{\mathcal U}_q(sl_2)$, the right integral $\mu$  is given by
\begin{equation}
\mu(E^i \, F^m \, K^n) = 
{\zeta}\, 
 \delta_{i, N-1} \, \delta_{m, N-1} \, \delta_{n, N+1}, 
 \label{eq:mu}
\end{equation}
where we choose the normalization factor as
\begin{equation}
\zeta =-\sqrt{\frac{2}{N}} \, ([N-1]!)^2
\label{eq:normalizationfactor}
\end{equation}
for future convenience.  
\medskip
\par\noindent
{\bf  Proposition 1.}
\begin{it}
The right integral satisfies
\begin{equation}
\mu(x \, y) = \mu(K^{1-N} \, y \, K^{N-1} \, x).  
\label{eq:rightintegral}
\end{equation}
\end{it}
\par\noindent
{\sc Proof.}
This comes \eqref{eq:mu} and the defining relations of $\overline{\mathcal U}_q(sl_2)$.  
\qed
\medskip
\par\noindent
{\bf Corollary 1.}  
\begin{it}
Let 
\begin{equation}
\phi(x) = \mu(K^{N+1} \, x),
\label{eq:fslf}
\end{equation} 
then $\phi(x \, y) = \phi(y\, x)$.  
\end{it}
\medskip\par\noindent
This $\phi$ is a fundamental tool for the construction of invariants of knots in three manifolds in this paper.  
\subsection{The universal $R$-matrix}
Let   $\mathcal A$ be the Hopf algebra generated by
$e$, $f$, $k$ and  relations
$$
\begin{aligned}
&k\, e\, k^{-1} = q \, e, \quad
k \, f \, k^{-1} = q^{-1} \, f, \quad
[e, f] = \dfrac{k^2 - k^{-2}}{q - q^{-1}},
\\&
e^N = f^N = 0, \quad k^{4N} = 1,
\quad
\epsilon(e) = \epsilon(f) = 0, \quad \epsilon(k) = 1, 
\\
&\Delta(e) = 1 \otimes e + e \otimes k^2, \quad
\Delta(f) = k^{-2} \otimes f + f \otimes 1, \quad
\Delta(k) = k \otimes k,
\\
&S(e) = -e \, k^{-2}, \quad S(f) = -k^2 \, f, \quad S(k) = k^{-1}.  
\end{aligned}
$$
Then there is an inclusion map $\iota : \overline{\mathcal U}_q(sl_2)\to \mathcal A$  given by 
\begin{equation}
\iota(E) = e, \qquad
\iota(F) = f, \qquad
\iota(K) = k^2.  
\label{eq:inclusion}
\end{equation}
In what follows, we often identify $E$ with $e$, $F$ with $f$, and $K$ with $k^2$.  
It is known that $\mathcal A$ is a ribbon quasitriangular Hopf algebra equipped with the universal $R$-matrix
\begin{equation}
\overline R = 
\dfrac{1}{4N}
\sum_{m=0}^{N-1}
\sum_{n,j=0}^{4N-1}
\dfrac{\{1\}^m}{[m]!} \,
q^{m(m-1)/2+m(n-j)-nj/2} \, e^m \, k^n \otimes f^m \, k^j.   
\label{eq:Rsmall}
\end{equation}
%
%
%
\subsection{Universal $\overline{\mathcal U}_q(sl_2)$ invariant}
\begin{figure}[htb]
$$
\raisebox{-5mm}{\begin{picture}(30,30)
\thicklines
\put(30,30){\vector(-1,-1){30}}
\put(0,30){\line(1, -1){12}}
\put(18, 12){\vector(1,-1){12}}
\put(8,8){\circle*{5}}
\put(22, 8){\circle*{5}}
\end{picture}}\quad
\longrightarrow\quad
 \sum_j a_j \otimes b_j,
\qquad
\raisebox{-5mm}{\begin{picture}(30,30)
\thicklines
\put(0,30){\vector(1,-1){30}}
\put(30,30){\line(-1, -1){12}}
\put(12, 12){\vector(-1,-1){12}}
\put(8,8){\circle*{5}}
\put(22, 8){\circle*{5}}
\end{picture}}
\quad\longrightarrow
\quad \sum_j b_j^\prime \otimes a_j^\prime,
$$
where $\overline R = \sum_j a_j \otimes b_j$ and $\overline R^{-1} = a_j^\prime \otimes b_j^\prime$,  
$$
\begin{matrix}
\quad
\raisebox{-2mm}{\epsfig{file=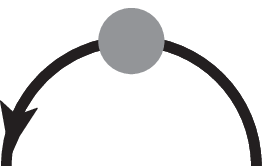, scale=0.3}} 
\to K^{-N+1},
\quad & 
\raisebox{-2mm}{\epsfig{file=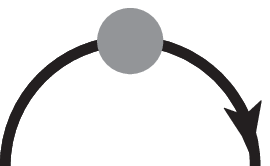, scale=0.3}}
\to 1,\quad & 
\raisebox{-2mm}{\epsfig{file=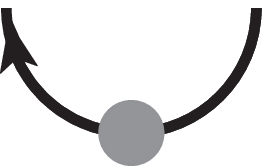, scale=0.3}}
\to K^{N-1},
\quad & 
\raisebox{-2mm}{\epsfig{file=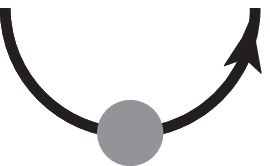, scale=0.3}}
\to 1.
\end{matrix}
$$
\caption{Universal invariant for crossings, maximal and minimal points}
\label{figure:universaldef}
\end{figure}

\begin{figure}[htb]
$$
L : \quad
\raisebox{-2.1cm}{\epsfig{file=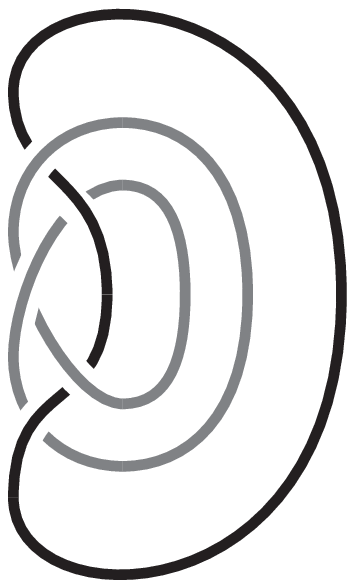, scale=0.8} }
\quad\longrightarrow\quad
\begin{matrix}
  \\{}\\
{}  \\
a_i \otimes b_i  \\
b_j^\prime \otimes a_j^\prime  \\[5pt]
a_k \otimes b_k  \\
b_l^\prime \otimes a_l^\prime  \\
a_m \otimes b_m  \\[-5pt]
&  \qquad x_2 \\{}\\
 &  \qquad x_1
\end{matrix}
\hspace{-13mm}
\raisebox{-2.1cm}{\epsfig{file=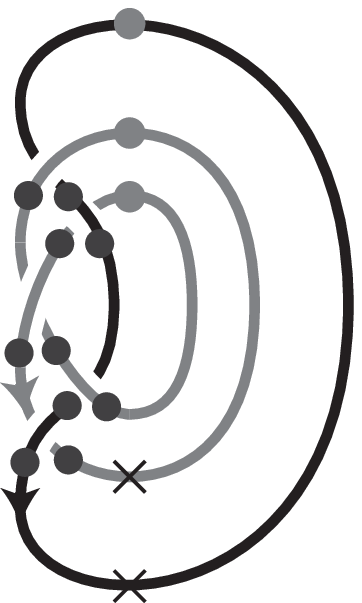, scale=0.8}}
\ %
\begin{matrix}
 \textcolor{black}{K^{1-N}} &  \\{}\\
\textcolor{black}{K^{1-N}} &  \\
\textcolor{black}{K^{1-N}} \\
{}\\{}\\{}\\{}\\{}\\{}\\{}\\
\end{matrix}
$$
$$
\longrightarrow \quad
\Psi_{x_1, x_2}(L) = \sum_{i, j, k, l, m}
a_m \, a_l^\prime \, a_j^\prime \, b_i \, K^{1-N} \otimes
b_m \, b_k \, b_j^\prime \, K^{1-N} \, a_l^\prime \, b_k \, a_i \, K^{1-N}
$$
\caption{Universal invariant for a link $L$}
\label{figure:universal}
\end{figure}

Let $L$ be diagram of a $k$-component framed oriented link $L= L_1\cup L_2 \cup \cdots L_r$ with blackboard framing given by a closed braid diagram. 
Assign the universal $R$ matrix or its inverse to each crossing and $K^{\pm(N-1)}$ or $1$ to each maximal and minimal points as in  Figure \ref{figure:universaldef}.   
Let $x_j$ be a point on $L_j$ other than the crossing points nor max/min points, and
we define $\Psi_{x_1, \cdots, x_r}(L)$  in ${\mathcal A}^{\otimes r}$
as follows.  
\begin{equation}
\Psi_{x_1, \cdots, x_r}(L)=
\sum_\nu u_1^{\nu}\otimes u_2^{\nu} \otimes\cdots\otimes u_r^\nu, 
\qquad
u_j^{\nu} = u_{j, 1}^{\nu} \, u_{j, 2}^{\nu} \, \cdots u_{j, p}^{\nu}
\quad \text{for $j=1$, $2$, $\cdots$, $r$}, 
\label{eq:universal}
\end{equation}
where $u_{j, 1}^{\nu}$,  $u_{j, 2}^{\nu}$, $\cdots$, $u_{j, p}^{\nu}$ are the elements
we meet when we walk through the component $L_j$ starting from $x_j$ to $x_j$ along its orientation as in Figure \ref{figure:universal}.  
For detail, see \cite{La} and  \cite{Oh2}.  
It is known that the element $\Psi_{x_1, \cdots, x_r}(L)$ is contained in the subspace $\iota(\overline{\mathcal U}_q(sl_2))^{\otimes r}$ of ${\mathcal A}^{\otimes r}$, and depends on the choices of $x_1$, $\cdots$, $x_r$.  
Let $\widehat{\overline{\mathcal U}}_q(sl_2)$ be the quotient 
$$
\widehat{\overline{\mathcal U}}_q(sl_2)
=
\overline{\mathcal U}_q(sl_2)/[\overline{\mathcal U}_q(sl_2), \overline{\mathcal U}_q(sl_2)],  
$$
and  $\psi(L)$ be the image of $\Psi_{x_1, \cdots, x_r}(L)$ in $\widehat{\overline{\mathcal U}}_q(sl_2)^{\otimes r}$, then the image $\psi(L)$ doesn't depend on the choices $x_1$, $\cdots$, $x_r$ and is an invariant of $L$.  
We  call $\psi(L)$ the {\it universal  $\overline{\mathcal U}_q(sl_2)$ invariant}.  
\subsection{Hennings invariants colored by symmetric linear functions and centers}
We recall Hennings' method in \cite{He} to retrieve numerical invariants from $\psi(L)$.  
\smallskip\par\noindent
{\bf Definition 1.}
An element $f$ in $\widehat{\overline{\mathcal U}}_q(sl_2)^*$ is called a {\it symmetric linear function} on $\overline{\mathcal U}_q(sl_2)$.  
In other words, $f$ is a linear functional on $\overline{\mathcal U}_q(sl_2)$ satisfying $f(xy) = f(yx)$. 
For example, the function $\phi$ on $\overline{\mathcal U}_q(sl_2)$ introduced in Corollary 1 can be considered as an element of $\widehat{\overline{\mathcal U}}_q(sl_2)^*$.  
\smallskip\par
For $f_1$, $f_2$, $\cdots$, $f_r \in \widehat{\overline{\mathcal U}}_q(sl_2)^*$,   
$$
(f_1\otimes f_2 \otimes \cdots \otimes f_r)(\Psi_{x_1, \cdots, x_r}(L))
=
\sum_\nu f_1(u_1^{\nu})\, f_2(u_2^{\nu}) \cdots f_r(u_r^\nu)
$$ 
depends only on $\psi(L)$ and is an invariant of $L$.  
Moreover, let $z_1$, $z_2$, $\cdots$, $z_r$ be elements of the center $\mathcal Z$ of $\overline{\mathcal U}_q(sl_2)$, then   $(f_1\otimes f_2 \otimes \cdots \otimes f_r)\big((z_1\otimes\cdots\otimes z_r)\Psi_{x_1, \cdots, x_r}(L)\big)$ is also an invariant of $L$, which we denote by $\psi_{(f_1, z_1), \cdots, (f_r, z_r)}(L)$.  
\par
\begin{figure}[htb]
 $$
 \begin{matrix}
 \begin{matrix}
 & \nearrow \hspace{-4mm}\swarrow& \raisebox{2mm}{$L \cup\,  \raisebox{-2mm}{\epsfig{file=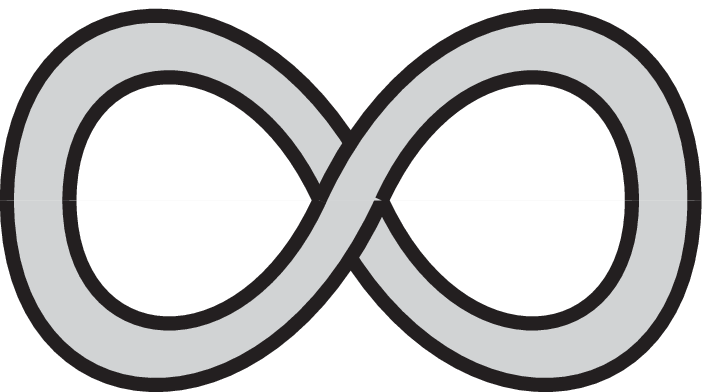, scale=0.2}}$\quad\ }
 \\
 L 
 \\
 & \searrow \hspace{-4mm} \nwarrow
 &  \raisebox{-2mm}{$L \cup\raisebox{-2mm}{ \epsfig{file=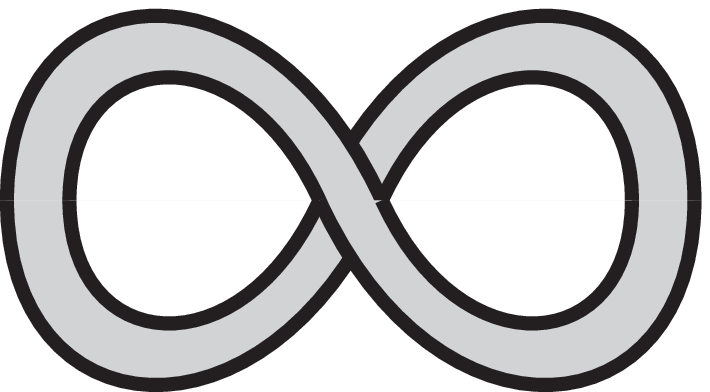, scale=0.2}}$\quad ,}
 \end{matrix}
& \qquad\quad &
 \raisebox{-8mm}{\epsfig{file=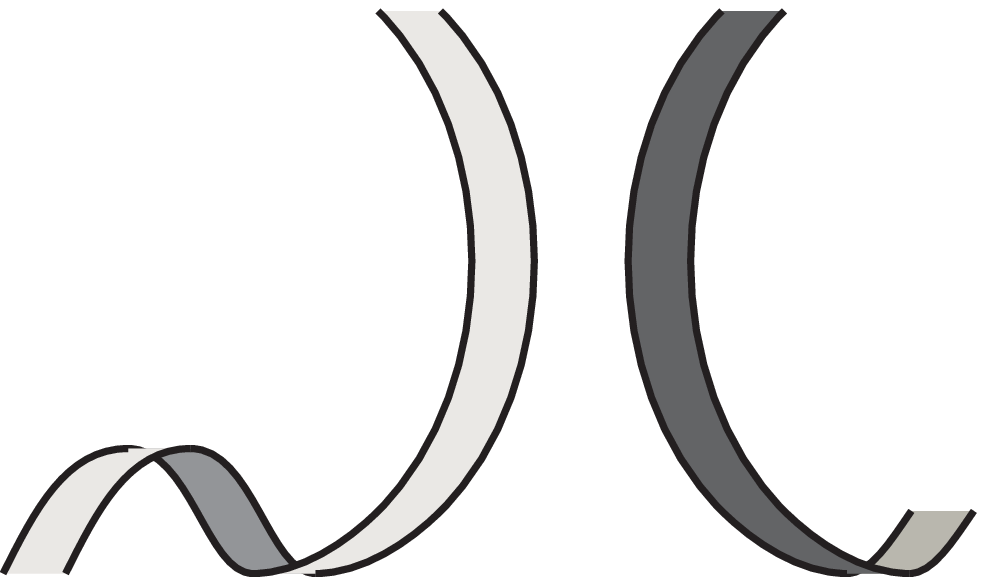, scale=0.3}}
 \longleftrightarrow
 \raisebox{-8mm}{\epsfig{file=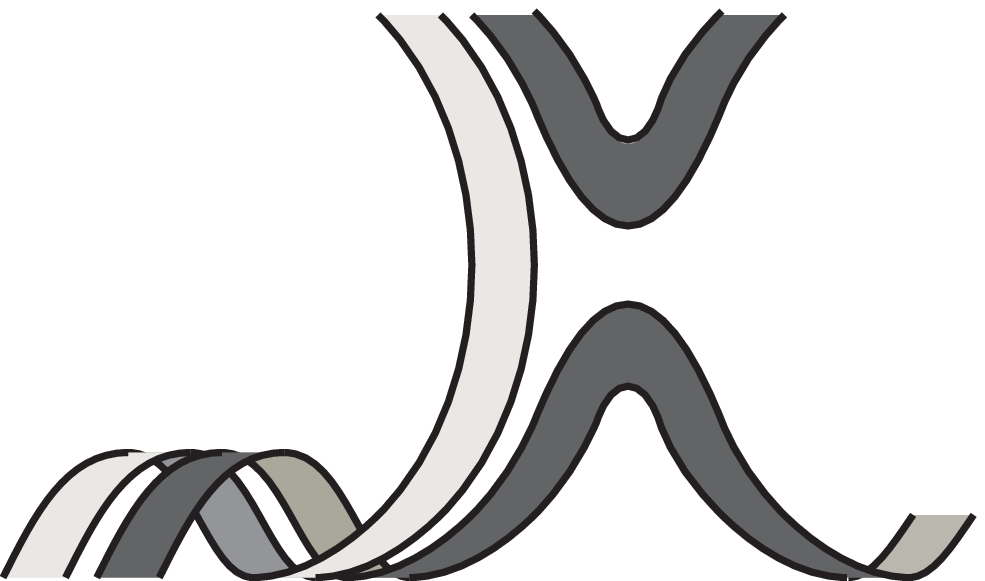, scale=0.3}}
 \\[30pt]
\text{${\mathcal O}_1$ move}
& & \text{${\mathcal O}_2$ move}
\end{matrix}
 $$
 \caption{${\mathcal O}_1$ and ${\mathcal O}_2$ moves}
 \label{fig:O2}
 \end{figure}

Hennings  shows that $\psi_{(\phi, 1), \cdots, (\phi, 1)}(L)$ is invariant under the second Kirby move ${\mathcal O}_2$ in Figure \ref{fig:O2},  
and we can construct a three manifold invariant from  $\psi_{(\phi, 1), \cdots, (\phi, 1r)}(L)$ by applying the normalization for the first Kirby move ${\mathcal O}_1$.  
Let $U_\pm$ be the unknot with $\pm 1$ framing.  
Let $s_+(L)$ (resp.  $s_-(L)$) be the number of positive (reps. negative) eigenvalues of the linking matrix of $L$.  
Here, the linking matrix $M = (m_{ij})_{1 \leq i,j \leq r}$ of $L$ is given by
$$
\begin{cases}
m_{ij} = \text{the linking number of $L_i$ and $L_j$}  \quad(i \neq j), \\
m_{ii} = \text{the writhe (the number indicating the framing) of $L_i$}. 
\end{cases}
$$
\smallskip
\par\noindent
{\bf Theorem 1.} (Hennings \cite{He})
\begin{it}
Let
 $$
H(M_L) = \dfrac{\psi_{(\phi, 1), \cdots, (\phi, 1)}(L)}
{\psi_{(\phi, 1)}(U_+)^{s_+(L)} \, \psi_{(\phi, 1)}(U_-)^{s_-(L)}}.  
$$
Then $H(M_L)$ is an invariant of the three manifold $M_L$ obtained from the surgery of $S^3$ along the framed link $L$.  
\end{it}
\smallskip
%
%
%
 %
%
%
\subsection{Colored Hennings invariants for links in three manifolds}
Let $M$ be an oriented three manifold given by the surgery along a framed link $L=L_1\cup L_2 \cup \cdots \cup L_p$ in $S^3$ and   
 $\widetilde K$ be a framed link in $M$ whose pre-image in $S^3$ is $\widehat K=K_1\cup K_2 \cup \cdots \cup K_r$ which does not intersect with $L$,   
where $L_i$ $(1 \leq i \leq p)$, $K_j$ $(1 \leq j \leq r)$ are the connected components of $L$ and $\widehat K$ respectively.  
For $z_1$, $\cdots$, $z_r$ in $\mathcal Z$ and symmetric linear functions $f_1$, $\cdots$, $f_r$ in $\widehat{\overline{\mathcal U}}_q(sl_2)^*$, 
we put
\begin{equation}
\psi_{(f_1, z_1), \cdots, (f_r, z_r)}(\widetilde K)
=
\dfrac{
\psi_{(f_1, z_1), \cdots, (f_r, z_r), (\phi, 1) \cdots, (\phi, 1)}(\widehat K\cup L)}
{\psi_{(\phi, 1)}(U_+)^{s_+(L)} \, \psi_{(\phi, 1)}(U_-)^{s_-(L)}}.
\label{eq:coloredHennings}
\end{equation}
\smallskip
\par\noindent
{\bf Theorem 2.}
\begin{it}
$\psi_{(f_1, z_1), \cdots, (f_r, z_r)}(\widetilde K)$ be an invariant of the link $\widetilde K$ in $M$ where the $i$-th component of $\widetilde K$ is colored by $(f_i, z_i)$ for $i = 1$, $2$, $\cdots$, $r$.  
\end{it}
\smallskip\par
\begin{figure}[htb]
$$
\begin{matrix}
L_j \raisebox{-4mm}{\epsfig{file=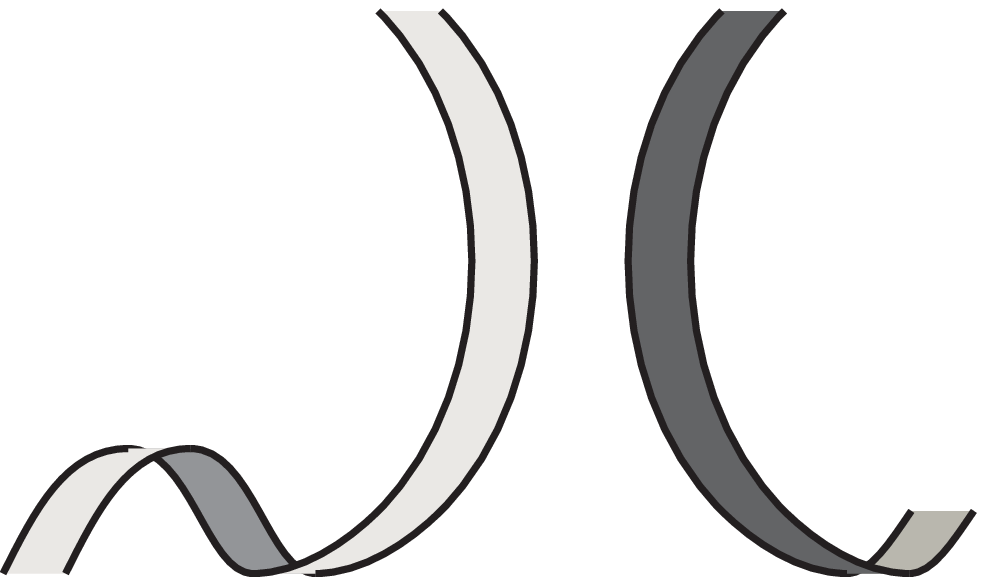, scale=0.3}}
K_i
\raisebox{4mm}{$ \longleftrightarrow$}
 L_j
 \raisebox{-4mm}{\epsfig{file=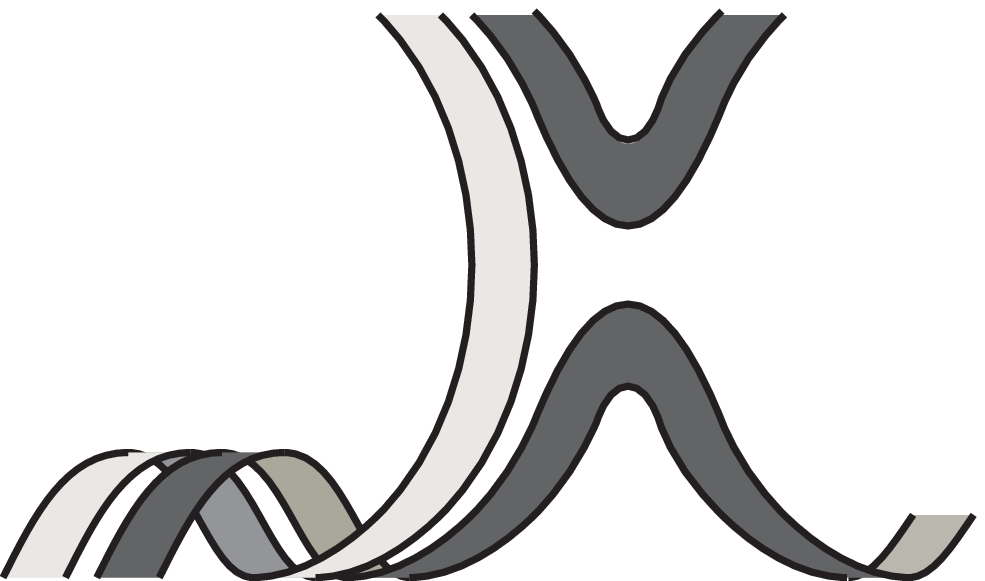, scale=0.3}}
 K_i
 \\[30pt]
\text{${\mathcal O}_2$ move for $\widehat K$ and $L$}
\end{matrix}
$$
\caption{Handle slide of a component $K_i$ of $\widehat K$ along a component $L_j$ of $L$}
\label{figure:pass}
\end{figure}

\noindent
{\sc Proof.}
We investigate the isotopy move of $\widetilde K$ by its pre-image $\widehat K$ in $S^3$.  
The isotropy of $\widetilde K$ which does not hit to the image of $L$ corresponds and isotopy of $\widehat K$ in $S^3$ which does not intersect with $L$.  
If a component $K_i$ of $\widehat K$ pass the image of a component $L_j$ of $L$ in $M$, then
the pre-image of this move in $S^3$ is given by the handle slide illustrated in Figure \ref{figure:pass}.  
Since $\phi$ is applied to all the components of $L$, 
$\psi_{(f_1, z_1), \cdots, (f_r, z_r), (\phi, 1) \cdots, (\phi, 1)}(\widehat K\cup L)$
does not change by this handle slide move.  
\qed
\section{Centers and symmetric linear functions of $\overline{\mathcal U}_q(sl_2)$}
In this section, we recall irreducible and indecomposable representations of $\overline{\mathcal U}_q(sl_2)$ and describe 
its centers and the symmetric linear functions explicitly.   
We also explain the $SL(2, {\mathcal Z})$ action on the centers of $\overline{\mathcal U}_q(sl_2)$.  
\subsection{Representations of  $\mathcal A$}
To explain representations of $\overline{\mathcal U}_q(sl_2)$, we first describe representations of the Hopf algebra ${\mathcal A}$ introduced in \S1.3.  
Let $U_s^{\alpha, \beta}$ be the $s$-dimensional irreducible representations of $\mathcal A$ labeled by $\alpha, \beta = \pm$ and $1 \leq s \leq N$.  
Let $t = \exp(\pi \, \sqrt{-1}/2N)$.  
The module  $U_s^{\pm,\pm}$ is  spanned by elements $u_n^{\pm,\pm}$ for $0 \leq n \leq s-1$, where the action of $\mathcal A$  is given by
$$
\begin{aligned}
k \, u_n^{\alpha, \beta} 
&= 
\beta\, \sqrt{\alpha}\, t^{s-1-2n} \, u_n^{\alpha, \beta}, 
\quad \text{$\sqrt{\alpha} = 1$ if $\alpha = +$ and $\sqrt{\alpha} = \sqrt{-1}$ if $\alpha = -$},
\\
e \, u_n^{\alpha, \beta} &= \alpha\, [n][s-n] \, u_{n-1}^{\alpha, \beta}, 
\quad  1 \leq n \leq s-1,
\qquad
e \, u_0^{\alpha, \beta} = 0, 
\\
f \, u_n^{\alpha, \beta} &= u_{n+1}^{\alpha, \beta}, 
\quad
0 \leq n \leq s-2,
\qquad\qquad\qquad\quad
f \, u_{s-1}^{\alpha, \beta} = 0 .  
\end{aligned}
$$
Especially, $U_1^{+,+}$ is the trivial module for which $k$ acts by 1 and  $e$, $f$ act by 0.  
The weights (eigenvalues of $k$) occurring in $U_s^{+,\pm}$ are
$$
\pm t^{s-1}, \ \pm t^{s-3}, \ \cdots , \ \pm t^{-s+1},
$$
and the  weights occurring in $U_{N-s}^{-, \pm}$ are
$$
\pm t^{2N-s-1}, \ \pm t^{2N-s-3}, \ \cdots, \ \pm t^{s+1}.  
$$
\par
Let  $V_s^{\alpha, \beta}$ $(1 \leq s \leq N)$ be the $N$ dimensional representation with highest-weight 
$\beta \, \sqrt{\alpha}\, t^{s-1}$ spanned by elements $v_n^{\pm,\pm}$ for $0 \leq n \leq N-1$, where the action of $\mathcal A$  is given by 
$$
\begin{aligned}
k \, v_n^{\alpha, \beta} 
&= 
\beta\, \sqrt{\alpha}\, t^{s-1-2n} \, v_n^{\alpha, \beta}, 
\quad \text{$\sqrt{\alpha} = 1$ if $\alpha = +$ and $\sqrt{\alpha} = \sqrt{-1}$ if $\alpha = -$},
\\
e \, v_n^{\alpha, \beta} &= \alpha\, [n][s-n] \, v_{n-1}^{\alpha, \beta}, 
\quad  1 \leq n \leq N-1,
\qquad
e \, v_0^{\alpha, \beta} = 0, 
\\
f \, v_n^{\alpha, \beta} &= v_{n+1}^{\alpha, \beta}, 
\quad
0 \leq n \leq N-2,
\qquad\qquad\qquad\quad
f \, v_{N-1}^{\alpha, \beta} = 0 .  
\end{aligned}
$$
Note that   
$V_N^\pm = U_N^\pm$.   
For $1 \leq s \leq N-1$,  $V_s^{\alpha, \beta}$ satisfies the exact sequence
$$
0 \longrightarrow  U_{N-s}^{-\alpha, -\beta} \longrightarrow
V_s^{\alpha,\beta} \longrightarrow U_s^{\alpha,\beta} \longrightarrow 0,   
$$
and there are projective modules 
$P_s^{\alpha,\beta}$ satisfying the following exact sequence.  
$$
0 \longrightarrow  V_{N-s}^{-\alpha,-\beta} \longrightarrow
P_s^{\alpha,\beta} \longrightarrow V_s^{\alpha,\beta} \longrightarrow 0.
$$
Actual description of the structure of  $\mathcal A$-modules $P_s^{\pm,\pm}$ is given  in \cite{JiMiTa}, 
which is based on the construction in \cite{ReTu}.    
The module $P_s^{+,\beta}$ $(\beta = \pm)$ has a basis
$$
\{x_j^{+,\beta},\ y_j^{+,\beta}\}_{0 \leq j \leq N-s-1}
\cup
\{a_n^{+,\beta}, \ b_n^{+,\beta}\}_{0 \leq n \leq s-1}.  
$$ 
The action of $k$ is given by
$$
\begin{aligned}
k \, x_j^{+,\beta} &= \beta\,t^{2N-s-1-2j} \, x_j^{+,\beta}, 
\ \ 
&k \, y_j^{+,\beta}=\beta\,t^{-s-1-2j}\, y_j^{+,\beta},
\qquad\qquad
&0 \leq j \leq N-s-1,
\\
k \, a_n^{+,\beta} &=\beta\, t^{s-1-2n} \, a_n^{+,\beta},
\qquad\quad
&k \, b_n^{+,\beta} = \beta\,t^{s-1-2n} \, b_n^{+,\beta}, 
\qquad\qquad\ \ 
&0 \leq n \leq s-1,
\end{aligned}
$$
The actions of $E$ and $F$ are given as follows.  
$$
\begin{aligned}
E \, x_j^{+,\beta}, &=
-[j][N-s-j] \, x_{j-1}^{+,\beta}, \quad 0 \leq j \leq N-s-1
\quad
(\text{with}\  x_{-1}^{+,\beta} =0),
\\
E\, y_j^{+,\beta} &= 
\begin{cases}
-[j][N-s-j] \, y_{j-1}^{+,\beta}, & 1 \leq k \leq N-s-1,
\\
a_{s-1}^{+,\beta}, & j=0,
\end{cases}
\\
E \, a_n^{+,\beta} &= 
[n][s-n] \, a_{n-1}^{+,\beta}, \qquad\qquad 0 \leq n \leq s-1
\quad\qquad \text{(with $a_{-1}^{+,\beta} = 0$)},
\\
E \, b_n^{+,\beta} &=
\begin{cases}
[n][s-n]\, b_{n-1}^{+,\beta} + a_{n-1}^{+,\beta}, & 1 \leq n \leq s-1,
\\
x_{N-s-1}^{+,\beta}, & n=0,
\end{cases}
\end{aligned}
$$
$$
\begin{aligned}
F\, x_j^{+,\beta} &= 
\begin{cases}
x_{j+1}^{+,\beta}, & 0 \leq j \leq N-s-2,
\\
a_0^{+,\beta}, & j = N - s - 1,
\end{cases}
\\
F\, y_j^{+,\beta} &=
y_{j+1}^{+,\beta}, \qquad 0 \leq j \leq N-s-2
\qquad
\text{(with $y_{N-s}^{+,\beta} = 0$)},
\\
F \, a_n^{+,\beta} &= a_{n+1}^{+,\beta},
\qquad 0 \leq n \leq s-1
\qquad\qquad
\text{(with $a_s^{+,\beta} = 0$)}, 
\qquad\qquad\quad
\\
F\, b_n^{+,\beta} &=
\begin{cases}
b_{n+1}^{+,\beta}, & 0 \leq n \leq s-2,
\\
y_0^{+,\beta}, & n = s-1.  
\end{cases}
\end{aligned}
$$
The $\mathcal A$-module $P_{N-s}^{-,\beta}$ is described as follows.  
$P_{N-s}^{-,\beta}$ has a basis
$$
\{x_j^{-,\beta},\ y_j^{-,\beta}\}_{0 \leq j \leq N-s-1}
\cup
\{a_n^{-,\beta}, \ b_n^{-,\beta}\}_{0 \leq n \leq s-1}.  
$$ 
The action of $\mathcal A$ is given by
$$
\begin{aligned}
k \, x_j^{-,\beta} &= \beta \, t^{-s-1-2j} \, x_j^{-,\beta}, 
\qquad
k \, y_j^{-,\beta} = \beta \, t^{-s-1-2j}\, y_j^{-,\beta},
\qquad
0 \leq j \leq N-s-1,
\\
k \, a_n^{-,\beta} &= \beta \, t^{s-1-2n} \, a_n^{-,\beta},
\qquad
K \, b_n^{-,\beta} = \beta \, t^{-2N+s-1-2n} \, b_n^{-,\beta}, 
\qquad
0 \leq n \leq s-1,
\end{aligned}
$$
$$
\begin{aligned}
E \, x_j^{-,\beta} &=
-[j][N-s-j] \, x_{j-1}^{-,\beta}, \quad 0 \leq k \leq N-s-1
\quad
(\text{with}\  x_{-1}^{-,\beta} =0),
\\
E\, y_j^{-,\beta} &= 
\begin{cases}
-[j][N-s-j] \, y_{j-1}^{-,\beta} + x_{j-1}^{-,\beta}, & 1 \leq j \leq N-s-1,
\\
a_{s-1}^{-,\beta}, & j=0,
\end{cases}
\\
E \, a_n^{-,\beta} &= 
[n][s-n] \, a_{n-1}^{-,\beta}, \qquad\qquad 0 \leq n \leq s-1
\quad\qquad \text{(with $a_{-1}^{-,\beta} = 0$)},
\\
E \, b_n^{-,\beta} &=
\begin{cases}
[n][s-n]\, b_{n-1}^{-,\beta}, & 1 \leq n \leq s-1,
\\
x_{N-s-1}^{-,\beta}, & n=0,
\end{cases}
\\
F\, x_j^{-, \beta} &=
x_{j+1}^{-,\beta}, \qquad 0 \leq j \leq N-s-2
\qquad
\text{(with $x_{N-s}^{-,\beta} = 0$)},
\\
F\, y_j^{-,\beta} &= 
\begin{cases}
y_{j+1}^{-,\beta}, & 0 \leq j \leq N-s-2,
\\
b_0^-, & j = N - s - 1,
\end{cases}
\qquad
F\, a_n^{-,\beta} =
\begin{cases}
a_{n+1}^{-,\beta}, & 0 \leq n \leq s-2,
\\
x_0^{-,\beta}, & n = s-1.  
\end{cases}
\\
F \, b_n^{-,\beta} &= b_{n+1}^{-,\beta},
\qquad 0 \leq n \leq s-1
\qquad\qquad
\text{(with $b_s^{-,\beta} = 0$)}.
\end{aligned}
$$
%
%
%
\subsection{Representations of $\overline{\mathcal U}_q(sl_2)$}
By composing the inclusion map $\iota: \overline{\mathcal U}_q(sl_2) \to {\mathcal A}$ given by \eqref{eq:inclusion} to the above representations of $\mathcal A$, we get representations of $\overline{\mathcal U}_q(sl_2)$.  
As representations of $\overline{\mathcal U}_q(sl_2)$, $X_s^{\pm, +}$ and $X_s^{\pm, -}$ are the same one for $X = U$, $V$ and $P$.  
Therefore, we write $U_s^\pm$, $V_s^\pm$ and $P_s^\pm$ for $U_s^{\pm, \beta}$, $U_s^{\pm, \beta}$ and $P_s^{\pm, \beta}$ respectively.  
\subsection{Symmetric linear functions}
It is shown in \cite{Ra} that there is a linear isomorphism between the center $\mathcal Z$ of $\overline{\mathcal U}_q(sl_2)$ and the space of symmetric linear functions $\widehat{\overline{\mathcal U}}_q(sl_2)^*$ given by
$
z \mapsto \mu(K^{N+1} \, z \, \bullet)
$, 
where $\mu$ is the right integral of   $\overline{\mathcal U}_q(sl_2)$.  
Hence, the dimension of $\widehat{\overline{\mathcal U}}_q(sl_2)^*$ is $3N-1$.  
A symmetric linear function which is not a trace of any semisimple representation is also called {\it pseudo-trace} in \cite {Mi}.  
The actual description of $\widehat{\overline{\mathcal U}}_q(sl_2)^*$ is given by Arike in \cite{Ar}, which is spanned by the following functions
 $T_0$, $T_N$, $T_1^\pm$, $\cdots$, $T_{N-1}^\pm$,  $G_1$, $\cdots$, $G_{N-1}$.  
 \begin{itemize}
\item
$T_0$ is the trace of the representation on $U_N^-$.  
\item
$T_N$ is the trace of the representation on $U_N^+$.  
\item
$T_s^\pm$ is the trace of the representation on  $U_s^\pm$ $(1 \leq s \leq N-1)$.  
\item
$G_s$ is the sum of the following two traces.  
One is the trace of the $s \times s$ submatrix of the representation matrix on $P_s^+$ at the block of the row for $a_n^+$ and columns for $b_m^+$ $(0 \leq n, m  \leq s-1)$,
 and the another one is the trace of  the $(N-s) \times (N-s)$ submatrix of the representation matrix on$P_{N-s}^-$ at rows of $x_k^-$ and columns of  $y_l^-$ $(0 \leq k, l \leq N-s-1)$.
\end{itemize}
The symmetric linear function $\phi$ introduced in Corollary 1 is explicitly written in \cite{Ar} as  follows.  
\smallskip
\par\noindent
{\bf Proposition 2.}
\begin{it}
The symmetric linear function $\phi$ is given by
\begin{equation}
\phi = 
\alpha_0 \, T_0 + \alpha_N \, T_N + 
\sum_{s=1}^{N-1}\left(
\alpha_s \, T_s + \beta_s \, G_s
\right),   
\label{eq.phi}
\end{equation}
where
$$
\begin{aligned}
T_s &= T_s^+ + T_s^-, 
\quad
\alpha_0 = -\dfrac{1}{N \sqrt{2N}},
\quad
\alpha_s =
\dfrac{(-1)^{s-1}\,
\{s\}_+}{N \sqrt{2N}
}
,
\\
\alpha_N &= \dfrac{(-1)^{N}}{N \sqrt{2N}},
\quad\ \
\beta_s =
\dfrac{(-1)^{s-1}\, [s]^2}{N \sqrt{2N}} . 
\end{aligned}
$$
\end{it}
\smallskip
\par\noindent
{\sc Proof.}
The coefficients $\alpha_0$, $\alpha_N$, $\beta_s$ $(1 \leq s \leq N-1)$ are obtained from those in \cite{Ar} by multiplying the normalization factor $\zeta$ in \eqref{eq:normalizationfactor}.  
The coefficient 
$\alpha_s^\pm$ $(1 \leq s \leq N-1)$ is given in \cite{Ar} as
$$
\alpha_s =
-\beta_s \, \left(
\sum_{l=1}^{s-1}\dfrac{1}{[l][s-l]} - 
\sum_{l=1}^{N-s-1} \dfrac{1}{[l][N-s-l]}
\right).
$$
A  computation shows that 
$$
\dfrac{1}{[l]\, [s-l]}
= 
[s]^{-1} \left({q^l}{[l]^{-1}} + {q^{l-s}}{[s-l]^{-1}}\right),
$$
which implies
$$
\sum_{l=1}^{s-1}{[l]^{-1}[s-l]^{-1}}
=
{[s]^{-1}}\sum_{l=1}^{s-1} 
{\{l\}_+}{[l]^{-1}}.
$$
Similarly, we have
$$
\sum_{l=1}^{N-s-1} {[l]^{-1}[N-s-l]^{-1}}
=
 -{[s]^{-1}}\sum_{l=s+1}^{N-1} 
{\{l\}_+}{[l]^{-1}}.
$$  
Since 
$
\sum_{l=1}^{N-1}
{\{l\}_+}{\{l\}^{-1}}
= 0 
$, 
we know that
$$
\sum_{l=1}^{s-1} 
{\{l\}_+}{[l]^{-1}}
+
\sum_{l=s+1}^{N-1} 
{\{l\}_+}{[l]^{-1}}
=
-
{\{s\}_+}{[s]^{-1}},
$$
which implies
$\alpha_s = 
(-1)^{s-1}\, 
\{s\}_+{N^{-1}\sqrt{2N}^{-1}}
$
.   \hfill\qed
\subsection{Centers of $\mathcal A$}
The center ${\mathcal Z}$ of $\overline{\mathcal U}_q(sl_2)$ is investigated in \cite{FeGaSeTi} and the center ${\mathcal Z}_{\mathcal A}$ of $\mathcal A$ is obtained similarly  as follows.  
\smallskip
\par\noindent
{\bf Proposition 3.}\ 
\begin{it}
The center ${\mathcal Z}_{\mathcal A}$ of  $\mathcal A$ is $5N-1$ dimensional.  
Its commutative algebra structure is described as follows.  
There are four special central idempotents $e_0^\pm$ and $e_N^\pm$, other central idempotents $e_s$, $1 \leq s \leq N-1$, and $4(N-1)$ elements $w_s^{\pm,\pm}$ $(1 \leq s \leq N-1)$ in the radical 
such that
$$
\begin{tabular}{ll}
\baselineskip20pt
$e_s^\alpha \, e_t^{\alpha'}= \delta_{s, t} \, \delta_{\alpha, \alpha'}\, e_s^\alpha,$
& $s, t = 0, 1, \cdots, N,\quad \alpha, \alpha' = \pm \text{ or empty}, $
\vspace{2mm}
\\
$e_s \, w_t^{\pm,\pm}= \delta_{s,t} \, w_t^{\pm,\pm}$,\qquad
& $0 \leq s \leq N, \ 1 \leq t \leq N-1,$
\vspace{2mm}
\\
$w_s^{\alpha,\beta} \, w_t^{\alpha',\beta'} = 0$,
& $1 \leq s,\ t \leq N-1$.  
\end{tabular}
$$
\end{it}
\par
The center $e_N^\pm$ acts on $U_N^{+, \pm}$ as an identity and acts as 0 on the other modules.  
$e_0^\pm$ acts on $U_N^{-, \pm}$ as identity and acts as 0 on the other modules.  
$e_s$ acts on $P_s^{+, +}$, $P_s^{+,-}$, $P_s^{-, +}$ and $P_{N- s}^{-, -}$ as identity and acts as 0 on the other modules.  
The center $w_s^{+, \pm}$ acts on $P_s^{+, \pm}$ by
$w_s^{+, \pm} \, b_n^{+, \pm}  = a_n^{+, \pm}$, $w_s^+ \, a_n^{+, \pm} = 0$, 
$w_s^{+, \pm} \, x_k^{+, \pm} = 0$, $w_s^{+, \pm} \, y_k^{+, \pm} = 0$, and
acts on the other modules as 0.  
Similarly, $w_s^{-,\pm}$ acts on $P_s^{-,\pm}$ by
$w_s^{-,\pm} \, y_k^{-,\pm}  = x_k^{-,\pm}$, $w_s^{-,\pm} \, x_k^{-,\pm} = 0$, 
$w_s^{-,\pm} \, a_n^{-,\pm} = 0$, $w_s^{-,\pm} \, b_n^{-,\pm} = 0$, and
acts on the other modules as 0.  
\par
The center $\mathcal Z$ of $\overline{\mathcal U}_q(sl_2)$ is spanned by
$\be_s$, $0 \leq s \leq N$, $\bw_s^\pm$, $1 \leq s \leq N-1$ whose images in $\mathcal A$ are
$$
\iota(\be_0) = e_0^+ + e_0^-, \quad
\iota(\be_N) = e_N^+ + e_N^-, \quad
\iota(\be_s) = e_s, \quad
\iota(\bw_s^\pm) = w_s^{\pm, +} + w_s^{\pm, -}.  
$$
Any central element $z$ in  ${\mathcal Z}$ is a linear combination of $e_s$, $w_s^\pm$ as follows.  
$$
z = \sum_{s=0}^N a_s(z) \, \be_s + 
\sum_{s=1}^{N-1}
\left(
b_s^+(z) \, \bw_s^+ + b_s^-(z) \, \bw_s^-
\right).  
$$
\section{Generalized logarithmic invariants}
Here we generalize the logarithmic invariant of a knot in $S^3$ to a knot in a three manifold.   
The logarithmic invariant is represented by the center corresponding to a $(1, 1)$-tangle of the knot, and we extend it by combining with the Hennings invariant. 
We show that there are some generalized logarithmic invariants which cannot be expressed by the colored Hennings invariant.  
\subsection{Center corresponding to a knot in a three manifold}
Let $M$ be a three manifold obtained by the surgery along a framed link $L = L_1 \cup \cdots \cup L_p$ and $\widetilde K$ be a knot or a link in $M$.  
Let $\widehat K=K_1 \cup \cdots\cup K_r$ be the pre-image of $\widetilde K$ in $S^3$ as the setting of Section 1.  
Let $T$ be the tangle obtained by cutting the first component $K_1$ of  $K_1 \cup \cdots\cup K_r \cup L_1 \cup \cdots \cup L_p$.  
The universal invariant of knots in \cite{La} and \cite{Oh1} is generalized to tangles in \cite{Oh2}, which we can apply to $T$. 
Let $x_2$, $\cdots$, $x_{r+p}$ be points on $K_2$, $\cdots$, $K_r$, $L_1$, $\cdots$, $L_p$ and 
let 
$$
\Psi_{x_2, \cdots, x_{p+r}}(T) = \sum_\nu u_1^\nu\otimes  \cdots \otimes u_r^\nu \otimes u_{r+1}^\nu \otimes \cdots \otimes u_{r+p}^\nu
$$ 
\begin{figure}[htb]
$$
T : \ 
\raisebox{-15mm}{\epsfig{file=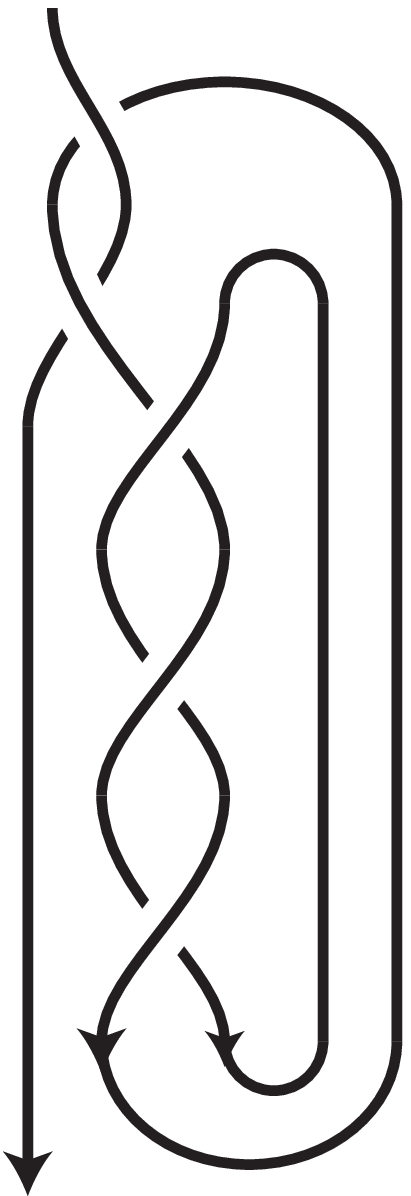, scale=0.3}}
\quad
\rightarrow
\ 
\begin{matrix}
b_i^\prime \otimes a_i^\prime \\
b_j^\prime \otimes a_j^\prime \\
a_k \otimes b_k \\[3pt]
a_l \otimes b_l \\[3pt]
a_m \otimes b_m\\{}
\end{matrix}
\raisebox{-2mm}{$\begin{matrix}\epsfig{file=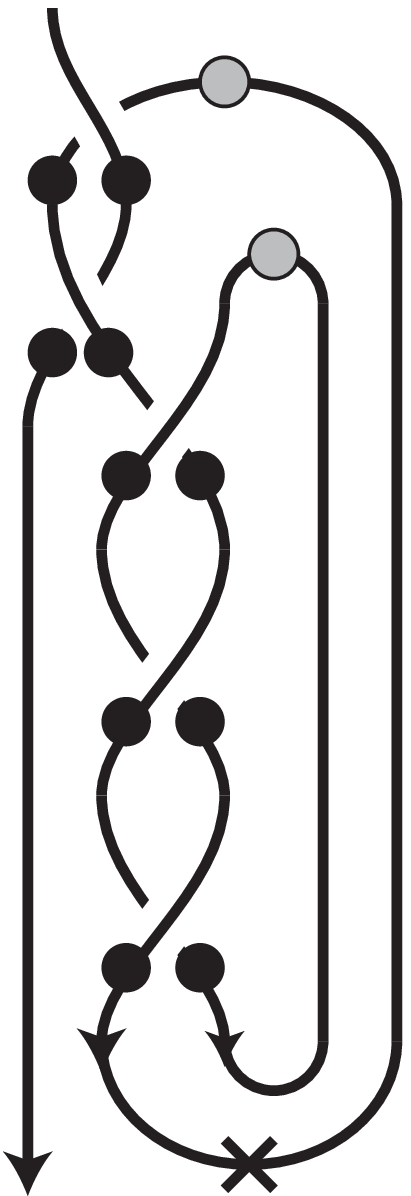, scale=0.3}\\[-5pt]
\quad x_2\end{matrix}$}\ \ 
\begin{matrix}
K^{1-N} \\
K^{1-N} \\[80pt]
\end{matrix}
\hspace{-8mm} \rightarrow\ 
\raisebox{-3mm}{\begin{tabular}{l}
$\Psi_{x_2}(T) = $\\[3pt]
\quad$\displaystyle\sum_{i,j,k,l,m}a_j^\prime \, b_i^\prime \otimes
a_m \, b_l \, a_k \, K^{1-N} \, b_m \, a_l \, b_k \, a_j^\prime \, b_i^\prime \, K^{1-N}$
\end{tabular}}
$$
\caption{Universal invariant for a tangle}
\end{figure}
be the element  of $\overline{\mathcal U}_q(sl_2)^{\otimes r+p}$ which is defined as \eqref{eq:universal}.  
For the component $K_1$, it is opened to make a tangle and we read the terms on this component from bottom to top.  
Let $(f_2, z_2)$, $\cdots$, $(f_r, z_r)$ be pairs of a symmetric linear function which are colors for the components $K_2$, $\cdots$, $K_r$ and 
 $s_+(L)$ (resp.  $s_-(L)$) be the number of positive (reps. negative) eigenvalues of the linking matrix of $L$.  
Then 
\begin{equation}
z_{(f_2, z_2), \cdots, (f_r, z_r)} (T)
= 
\psi_{(\phi, 1)}(U_+)^{-s_+(L)} \, \psi_{(\phi, 1)}(U_-)^{-s_-(L)}\,\sum_\nu \left(\prod_{i=2}^r f_i(z_i \, u_{i}^\nu)\prod_{j=1}^p \phi(u_{r+j}^\nu) \right) \, u_1^\nu
\label{eq:z}
\end{equation}
is contained in the center $\mathcal Z$ of $\overline{\mathcal U}_q(sl_2)$.  
\smallskip\par\noindent
{\bf Theorem 3.}
\begin{it}
The center
$z_{(f_2, z_2), \cdots, (f_r, z_r)} (T)$  
is an invariant of the colored knot $\widetilde K$ with a specified component $K_1$.  
\end{it}
\smallskip\par\noindent
{\bf Definition 2.}
Fix a basis $\{c_1, c_2, \cdots\}$ of the center  $\mathcal Z$.
Then the center
$z_{(f_2, z_2), \cdots, (f_r, z_r)} (T)$  
is expressed as a linear combination of this basis, and the coefficients are
also invariants of $\widetilde K$.  
We call these coefficients the {\it logarithmic invariants} of $\widetilde K$ since they are related to the logarithmic TQFT constructed in \cite{FeGaSeTi}.  
\medskip\par
The colored Hennings invariant $\psi_{(f_1, z_1), \cdots, (f_r, z_r)} (\widetilde K)$
is expressed as 
$$
\psi_{(f_1, z_1), \cdots, (f_r, z_r)} (\widetilde K)
=
f_1(K^{N+1}\, z_1\, z_{(f_2, z_2), \cdots, (f_r, z_r)} (T)).  
$$
Therefore, we have the following.  
\smallskip\par\noindent
{\bf Corollary 2.}
\begin{it}
The colored Hennings invariant with colors $(f_1, z_1)$, $\cdots$, $(f_r, z_r)$ is determined by the logarithmic invariants.  
\end{it} 
\subsection{Simplification of the coloring}
For the coloring of each component of a link, we use a pair $(f, z)$ where $f$ is a symmetric linear function and $z$ is a center.  
However, we know that 
$G_s(\bw_s^\pm\, u) =  T_s^\pm(\be_s\, u) = T_s^\pm(u)$, 
$T_s(\bw_s^\pm,\, u) = 0$, 
$G_s(\be_s\, u) = G_s(u)$, 
and so the coloring by $(G_s, \bw_s^\pm)$, $(T_s^\pm, \be_s)$ are equal to the coloring by $(T_s^\pm, 1)$,  the coloring by $(G_s, \be_s)$ is equal to the coloring by $(G_s, 1)$, and the coloring by $(T_s^\pm, \bw_s^\pm)$ vanishes.  
This means that the invariant with any coloring can be expressed as a linear combination of invariants with colorings $(T_s^\pm, 1)$ and  $(G_s, 1)$.  
Therefore, from now on, we use the coloring  by symmetric linear functions only.  
For example, $\psi_{f_1, \cdots, f_r}$ means $\psi_{(f_1, 1), \cdots, {(f_r, 1)}}$.  
\subsection{Coefficients of the basis}
Let $\widetilde K$ be an $r$ component framed link in a three manifold $M$ given by the surgery along a $p$ component framed link $L=L_1 \cup \cdots \cup L_p$, $\widehat K =K_1\cup \cdots \cup K_r$  be the pre-image of $\widetilde K$, and $T$ be a $(1,1)$-tangle obtained from $K\cup L$ by cutting the component $K_1$ as before.
For such $\widehat K$,  we have constructed the center $z_{f_2, \cdots, f_r}(T)$ (with simplified colorings) in Theorem 3, which is an invariant of $\widetilde K$ with specialized component $\widetilde K_1$.  
This element is expressed as a linear combination of the basis of $\mathcal Z$ as
\begin{multline}
z_{f_2, \cdots, f_r}(T) =
 a_{o, f_2, \cdots, f_r}(T) \, \be_0 + a_{N, f_2, \cdots, f_r}(T) \, \be_N + 
\\
\sum_{s=1}^{N-1}
\left(
a_{s, f_2, \cdots, f_r}(T) \, \be_s + b_{s, f_2, \cdots, f_r}^+(T) \, \bw_s^+ + b_{s, f_2, \cdots, f_r}^-(T) \, \bw_s^-
\right).  
\label{eq:coeff}
\end{multline}
From Theorem 3 and the definition of the logarithmic invariants, we have the following.  
\smallskip
\par\noindent
{\bf Corollary 3.}
\begin{it}
The coefficients $a_{s, f_2, \cdots, f_r}(T)$ $(0 \leq s \leq N)$, $b_{s, f_2, \cdots, f_r}^\pm(T)$ $(1 \leq s \leq N-1)$ are logarithmic invariants of $\widetilde K$ in $M$.  
\end{it}
\medskip\par
Now, let us compare the colored Hennings invariants 
$\psi_{f_1, \cdots, f_r}(\widetilde K)$ 
and the logarithmic invariants
$a_{s, f_2, \cdots, f_r}(T)$ $(0 \leq s \leq N)$, $b_{s, f_2, \cdots, f_r}^\pm(T)$ $(1 \leq s \leq N-1)$ 
coming from 
$z_{f_2, \cdots, f_r}(T)$.  
Since 
$$
T_0(K^{N+1}\, {\be_0})=T_N(K^{N+1}\, {\be_N}) = 0, 
$$ 
$$
T_s^\pm(K^{N+1}\, \be_s) = \pm[s], \quad
G_s(K^{N+1} \, \bw_s^\pm) = \mp(-1)^{s} \, [s]
\quad
\text{ for $1 \leq s \leq N-1$},
$$ 
we have
$$
\begin{aligned}
\psi_{T_0, f_2, \cdots, f_r}(\widetilde K) &=
\psi_{T_N, f_2, \cdots, f_r}(\widetilde K) =
\psi_{T_s^\pm, \,f_2, \cdots, f_r}(\widetilde K) = 0,
\\
\psi_{T_s^\pm,\, f_2, \cdots, f_r}(\widetilde K) &= \pm[s] \, a_{s, f_2, \cdots, f_r}(T), 
\\
\psi_{G_s, f_2, \cdots, f_r}(\widetilde K) &=
 (-1)^{s+1} \, [s] \, (b_{s, f_2, \cdots, f_r}^+(T)- b_{s, f_2, \cdots, f_r}^-(T)), 
\\
\psi_{G_s, f_2, \cdots, f_r}(\widetilde K) &=
\mp(-1)^{s} \, [s] \, a_{s, f_2, \cdots, f_r}(T).  
\end{aligned}
$$
These relations imply that the colored Hennings invariants are linear combinations of $a_{s, f_2, \cdots, f_r)}(T)$ and $b_{s, f_2, \cdots, f_r}^+(T)- b_{s, f_2, \cdots, f_r}^-(T)$ for $s = 1$, $\cdots$, $N-1$.  
However, the above relations don't determine
the invariants 
$a_{0, f_2, \cdots, f_r}(T)$, 
$a_{N, f_2, \cdots, f_r}(T)$, 
$b_{s, f_2, \cdots, f_r}^\pm(T)$ of $\widetilde K$ from the colors Hennings invariants.  
%
%
%
\section{Generalized Kashaev invariant}
In this section, we show that certain logarithmic invariants for links in a three manifold are generalizations of Kashaev's invariants for knots in $S^3$.  
\subsection{Colored Alexander invariants}
The colored Alexander invariant introduced in \cite{ADO} can be  constructed  from the quantum $R$-matrix of the medium quantum group $\widetilde{\mathcal U}_q(sl_2)$ as in \cite{JM}.  
The medium quantum group is defined as follows. 
\begin{multline*}
\widetilde{\mathcal U}_q(sl_2)
=
\left<
K,\ K^{-1},\ E,\ F \mid
K \, E \, K^{-1} = q^2 \, E,\ 
K \, F \, K^{-1} =q^{-2} \, F,\right.\\ 
[E, F] = \dfrac{K-K^{-1}}{q-q^{-1}}, \  
\left.
\vphantom{\dfrac{K-K^{-1}}{q-q^{-1}}}
E^N = F^N = 0
\right>.  
\end{multline*}
To see the relation between the colored Alexander invariant and the logarithmic invariants, we check the correspondence of the quantum $R$-matrices for $\overline{\mathcal U}_q(sl_2)$ and $\widetilde{\mathcal U}_q(sl_2)$.  
Let $\widetilde R$ be the $R$-matrix for the colored Alexander invariant, then it is given in \cite{JM}  by
\begin{equation}
\widetilde R
=
q^{\frac{1}{2}H \otimes H} \, 
\sum_{m=0}^{N-1} \dfrac{\{1\}^{m}}{\{m\}!} \, 
q^{\frac{m(m-1)}{2}} \, \left(E^m \otimes F^m\right),   
\label{eq:Rmiddle}
\end{equation}
where $H$ is a formal element satisfying $q^H = K$.  
\smallskip\par\noindent
{\bf Definition 3.}
A representation $V$ of $\widetilde{\mathcal U}_q(sl_2)$ is called an {\it integral weight representation} iff there are integers $\lambda_1$, $\cdots$, $\lambda_r$ such that $V = \oplus_i V_{\lambda_i}$ where $K\, v_i = q^{\lambda_i}$ for any $v_i \in V_{\lambda_i}$.  
\smallskip\par\noindent
{\bf Lemma 2.}
\begin{it}
Let $U$ and $V$ be integral weight representations, then the 
representations of $\widetilde R$ in \eqref{eq:Rmiddle} and $\overline R$ in \eqref{eq:Rsmall} on $U\otimes V$ are equal.  
\end{it}
\smallskip\par\noindent
{\sc Proof.}
Let $u\in U$ and $v\in V$ be weight vectors such that $K \, u = q^r \, u$ and $K \, v = q^s$.  
Then $q^{\frac{1}{2}H\otimes H}\, u \otimes v = q^{\frac{1}{2} rs} \, u \otimes v$,  
$$
\begin{aligned}
\widetilde R(u\otimes v)
&=
\sum_{m=0}^{N-1} \dfrac{\{1\}^m}{[m]!} \, q^{m(m-1)/2}\, q^{\frac{1}{2}H \otimes H}
E^m \, u \otimes F^m \, v
\\&=
\sum_{m=0}^{N-1} \dfrac{\{1\}^m}{[m]!} \, q^{m(m-1)/2 + (r+2m)(s-2m)/2}\,
E^m \, u \otimes F^m \, v,
\end{aligned}
$$
and
$$
\begin{aligned}
\overline R \,(u\otimes v) 
&= 
\dfrac{1}{4N} \sum_{m=0}^{N-1}
\sum_{n, j=0}^{4N-1}
\dfrac{\{1\}^m}{[m]!} \, 
q^{m(m-1)/2 + (r/2+m)n + (s/2-m-n/2)j} \, 
E^m \, u \otimes F^m\, v
\\&=
\sum_{m=0}^{N-1}
\dfrac{\{1\}^m}{[m]!} \, 
q^{m(m-1)/2 + (r+2m)(s-2m)/2} \, 
E^m \, u \otimes F^m\, v.
\end{aligned}
$$
Hence the actions of the $R$-matrices are equal.  
\qed
\smallskip\par
The colored Alexander invariant is defined for non-integral representations, but 
this lemma shows that $\widetilde R$ is also well-defined for integral weight representations.  
Let $\widetilde K$ be a knot in a three manifold $M$ given by the surgery along a framed link $L$ in $S^3$, $\widehat K$ be the pre-image of $\widetilde K$ in $S^3$,  and $T$ be a tangle obtained from $\widehat K\cup L$ by cutting a component $K_1$ of $\widehat K$.  
Let $\sum_\nu u_1^\nu\otimes \cdots \otimes u_{r+p}^\nu$ be the universal invariant of the tangle $T$ in $\left(\widetilde{\mathcal U}_q(sl_2)/\left[\widetilde{\mathcal U}_q(sl_2), \widetilde{\mathcal U}_q(sl_2)\right]\right)^{\otimes(r+p)}$ constructed from $\widetilde R$.  
In this case we  use  the element $H$ and  infinite sums which converge on any finite dimensional representations.  
\smallskip\par\noindent
{\bf Remark 2.}
Another construction of a universal invariant corresponding to the colored Alexander invariant is given by Ohtsuki in \cite{Oh1} by using colored ribbon Hopf algebras.  
\smallskip
\par
For $\lambda\in {\mathbb C}$, let $\rho^\lambda$ be the highest weight representation of $\widetilde{\mathcal U}_q(sl_2)$ with highest weight $\lambda-1$ and $\chi^{\lambda}$ be the character of $\rho^{\lambda}$, i.e. the trace of the representation matrix of $\rho^{\lambda}$.  
Let ${\mathcal X}(\lambda)$ be the representation space of $\rho^\lambda$ spanned by the weight vectors $v_0^\lambda$, $v_1^\lambda$, $\cdots$, $v_{N-1}^\lambda$, on which
$\widetilde{\mathcal U}_q(sl_2)$ acts by
$$
K \, v_n^\lambda = 
q^{\lambda-1-2n}\, v_n,
\qquad
E \, v_n^\lambda = 
[n][\lambda-n]\, v_{n-1},
\qquad
F \, v_n = 
\, v_{n+1},
$$
where $\, v_N^\lambda=0$.  
Then ${\mathcal X}(\lambda)$ is irreducible if $\lambda \in ({\mathbb C}\setminus {\mathbb Z}) \cup N{\mathbb Z}$.  
For $\lambda_1$, $\cdots$, $\lambda_r \in {\mathbb C}$, 
$
\sum_\nu
\rho_{\lambda_1}(u_1^\nu)\left(\prod_{i=2}^{r}\chi_i(u_i^\nu)\prod_{j=1}^p\phi(u_{r+j}^\nu)\right)  
$
is a scalar matrix and let $A_{\lambda_1, \cdots, \lambda_r, \phi, \cdots, \phi}(T)$ be the corresponding scalar.  
\smallskip\par\noindent
{\bf Remark 3.} 
If the weight $\lambda_1$, $\cdots$, $\lambda_r$ are all specialized to integers, then $A_{\lambda_1, \cdots, \lambda_r, \phi, \cdots, \phi}(T)$ coincides with the logarithmic invariant $a_{s_1, T_{s_2}, \cdots, T_{s_r}}(T)$ defined by \eqref{eq:coeff},  where $s_i \equiv \lambda_i$ or $2N-\lambda_i \ ({\rm mod} \ 2N)$.  
\smallskip
\par
If $L$ is empty and $\widehat  K$ is a framed link in $S^3$,  we know the following for the tangle $T$ corresponding to $\widehat K$ obtained by cutting the  component $K_1$ of $\widehat K$.  
\smallskip\par\noindent
{\bf Theorem 4\ }\cite{ADO}, \cite{JM}.
\begin{it}
For $\lambda_1$, $\cdots$, $\lambda_r$ in ${\mathbf C}\setminus {\mathbf Z}$, 
let
$$
\ADO_{\lambda_1, \cdots, \lambda_r}(T)
=
\dfrac{\sin(\lambda_1 \, \pi/N)}
{\sqrt{-1}^{N-1} \, \sin \lambda_1\pi}\, 
 A_{\lambda_1, \cdots, \lambda_r}(T).  
$$
Then $\ADO_{\lambda_1, \cdots, \lambda_r}(T)$ is an invariant of the link $\widehat K$ in $S^3$ which does not depend on the component $K_1$.  
\end{it}
\smallskip\par\noindent
{\bf Remark 4.}
For a framed link $\widehat K$ in $S^3$, $A_{N, \cdots, N}(T)$ is equal to Kashaev's invariant, which is equal to the colored Jones invariant corresponding to the $N$ dimensional representation of ${\mathcal U}_q(sl_2)$ at $q =\exp(\pi \sqrt{-1}/N)$.  
\subsection{Generalized Kashaev invariant}
We introduce generalized Kashaev invariants for links in three manifolds as a generalization of $A_{N, \cdots, N}(T)$ combining with the Hennings invariant.  
\smallskip\par\noindent
{\bf Theorem 5.}
\begin{it}
Let $\widetilde K$ be a link in a three manifold $M$ given by the surgery along a framed link $L$ and $\widehat K$ be the pre-image of $\widehat K$ in $S^3$.  
Let $T$ be the tangle obtained from $\widehat K \cup L$, and $s_+(L)$, $s_-(L)$ be the numbers of positive and negative eigenvalues of the linking matrix of $L$.  
Then  
$$
\psi_\phi(U_+)^{-s_+(L)}\, \psi_\phi(U_-)^{-s_-(L)} \, A_{\underset{r}{\underbrace{\scriptstyle N, \cdots, N}}, \ \underset{p}{\underbrace{\scriptstyle \phi, \cdots, \phi}}}(T)
$$
does not depend on the choice of the specified component $K_1$ of $\widetilde K$ to make the tangle $T$, and is an invariant of $\widetilde K$.  
\end{it}
\medskip\par
\begin{figure}[htb]
$$
\epsfig{file=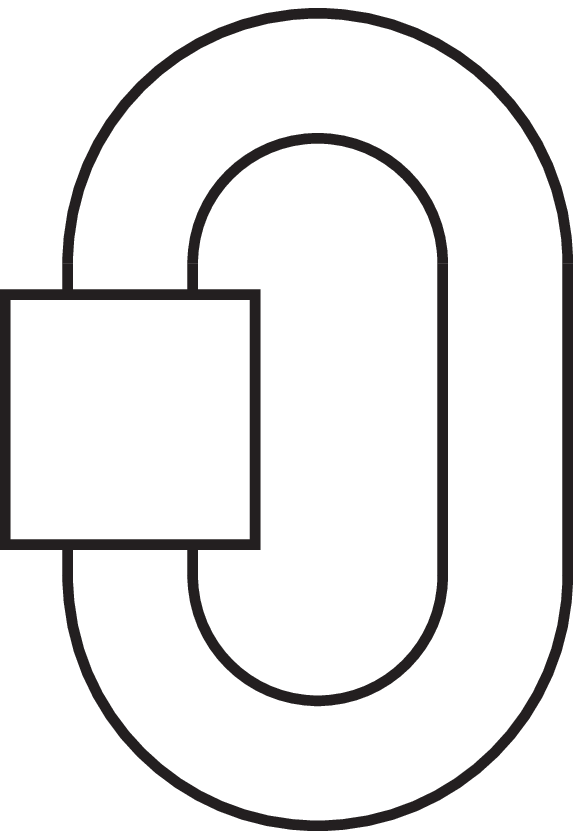, scale=0.4}
\hspace{-20mm}\raisebox{15mm}{$T_2$}\raisebox{20mm}{\quad\,$K_1$\quad\ \,$K_2$}
\qquad\qquad
\raisebox{8mm}{$\raisebox{8mm}{$\rho : \quad$}
\epsfig{file=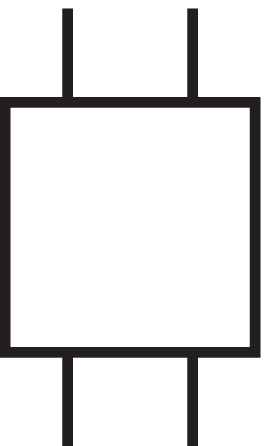, scale=0.4}
\hspace{-7mm}\raisebox{8mm}{$T_2$}\hspace{-12mm}\raisebox{20mm}{$V_N^{(+,+)}$\ \,$V_N^{(+,+)}$}
\raisebox{8mm}{$\longrightarrow\quad \rho(T^{(2)}) \in \End(V_N^{(+,+)}\otimes V_N^{(+,+)})$}$}
$$
\caption{The tangle $T^{(2)}$}
\label{figure:tangle}
\end{figure}

\begin{figure}[htb]
$$
\hspace{-15mm}
\raisebox{14mm}{$T_1^{(2)}$ :\ }  \quad
\epsfig{file=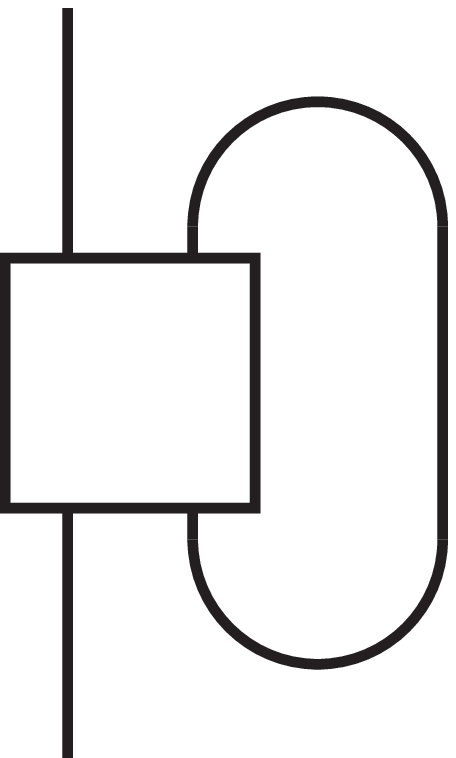, scale=0.4}
\hspace{-16mm}\raisebox{14mm}{$T^{(2)}$}
\qquad\ \raisebox{20mm}{$K_2$}
\hspace{-28mm}\raisebox{30mm}{$K_1$}
\qquad\qquad\qquad\qquad\qquad
\raisebox{14mm}{$T_2^{(2)}$ :\  }  \quad
\epsfig{file=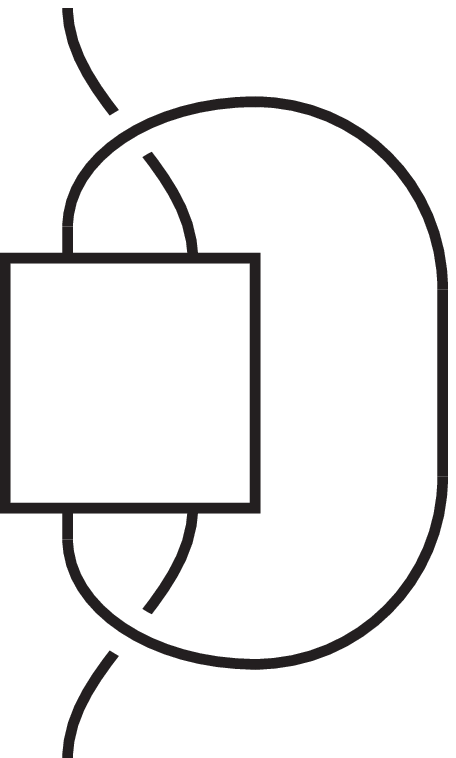, scale=0.4}
\hspace{-17mm}\raisebox{14mm}{$T^{(2)}$}
\qquad\ \raisebox{20mm}{$K_1$}
\hspace{-28mm}\raisebox{30mm}{$K_2$}
\hspace{-5mm}\raisebox{14mm}{$\begin{matrix}
\sigma \\[17mm] \sigma^{-1}
\end{matrix}$}
$$
\caption{The $(1,1)$-tangles $T_1^{(2)}$ and $T_2^{(2)}$}
\label{figure:twotangle}
\end{figure}

\noindent
{\sc Proof.}
We show that $A_{\underset{r}{\underbrace{\scriptstyle N, \cdots, N}}, \ \underset{p}{\underbrace{\scriptstyle \phi, \cdots, \phi}}}(T)$ does not depend on the choice of the component $K_1$.  
We assume that the number of the components of $\widehat K$  is greater than one.  
Let $T^{(2)}$ be a $(2,2)$-tangle obtained from $\widehat K$ as in Figure \ref{figure:tangle}.  
We associate the representation $V_N^{+,+}$ to the components of $\widehat K$.  
Then $T^{(2)}$ corresponds to an element $\rho(T^{(2)}) \in \End_{{\mathcal U}_q(sl_2)}(V_N^{+,+} \otimes V_N^{+,+})$, where $V_N^{+,+} \otimes V_N^{+,+}$ is split into a direct sum of indecomposable   $\overline{\mathcal U}_q(sl_2)$ modules as follows. 
\begin{equation}
V_N^{+,+} \otimes V_N^{+,+}
=
\begin{cases}
\displaystyle\bigoplus_{s=0}^{(N-2)/2} {\mathcal P}_s^+, & \text{if $N$ is even,}
\\
V_N^{+,+} \oplus \displaystyle\bigoplus_{s = 0} ^{(N-3)/2} {\mathcal P}_s^+, & \text{if $N$ is odd.}
\end{cases}
\label{eq:split}
\end{equation}
Note that this  is a multiplicity-free decomposition.  
Hence the action of $\rho(T^{(2)})$ is decomposed into a direct sum of the actions on   ${\mathcal P}_s$ and $V_N^{+,+}$ which commute with the action of $\overline{\mathcal U}_q(sl_2)$.  
Let $T_1^{(2)}$ and $T_2^{(2)}$ be $(1,1)$-tangles obtained by closing the right strings of $T^{(2)}$ and $\sigma \, T^{(2)} \, \sigma^{-1}$ as in Figure \ref{figure:twotangle}.  
Now compare the scalar corresponding to $T_1^{(2)}$ and $T_2^{(2)}$.  
The action of $\sigma$ on $V_N^{+,+} \otimes V_N^{+,+}$ commutes with the action of $\overline{\mathcal U}_q(sl_2)$.   
Therefore, the images  of $\sigma$ and $T^{(2)}$  in $\End(P_s^{+,+})$ and $\End(V_n^{+,+})$ are both contained in the commutants  with respect to $\overline{\mathcal U}_q(sl_2)$.
From the construction of representations of $\overline{\mathcal U}_q(sl_2)$, 
it is easy to see  that the above commutants are all abelian.  
This implies that   
 the images of $T_1^{(2)}$,  $T_2^{(2)}$  in $\End({\mathcal P}_s^{+,+})$ and 
$\End(V_N^{+,+})$ are the same ones, and then we get 
$$
A_{N, \cdots, N, \phi, \cdots, \phi}(T_1^{(2)})
=
A_{N, \cdots, N, \phi, \cdots, \phi}(T_2^{(2)}).  
$$
\qed
\par
By using $A_{N, \cdots, N, \phi, \cdots, \phi}$, we introduce the generalized Kashaev invariant ask follows.  
\smallskip\par\noindent
{\bf Definition 4.}
Let 
$$
\widetilde \GK_N(T) = A_{N, \cdots, N, \phi, \cdots, \phi}(T),
$$ 
which is equal to 
$a_{N, T_N, \cdots, T_N, \phi, \cdots, \phi}(T)$ by Remark 3, 
 and let
$$
\GK_N(\widetilde K) = \dfrac{\widetilde \GK_N(T)}{\psi_\phi(U_+)^{-s_+(L)}\, \psi_\phi(U_-)^{-s_-(L)}}.
$$  
We call $
\GK_N(T)
$
the {\it generalized Kashaev invariant} of $\widetilde K$.  
\smallskip\par
To get more computable expression of $\GK_N(\widetilde K)$, we express the symmetric linear function $G_s$ by  derivatives of the diagonal elements of $\rho^\lambda$.  
%
%
%
%
%
%
\subsection{An expression of $G_s$}
We first introduce a non-irreducible module of the medium quantum group $\widetilde{\mathcal U}_q(sl_2)$ which is isomorphic to direct sum of two non-integral highest weight representations.  
Let $t$ be an integer with $1 \leq s \leq p$ and ${\mathcal Y}(\lambda, s)$ be the $\sresq$ module which is spanned by 
weight vectors $c_j$ and $d_j$ for $0 \leq j \leq N-1$.
The action of $\sresq$ is given by
$$
\begin{aligned}
K \, c_n
&=
q^{\lambda-1-2n} \, c_n, 
\qquad
K \, d_n
=
q^{\lambda-1-2s-2n} \, d_n, 
\qquad
0 \leq n \leq N-1,
\end{aligned}
$$
$$
\begin{aligned}
E\, c_n
&=
\begin{cases}
0, & \quad n = 0, \\
[n][\lambda-n]\,c_{n-1} , 
&
\quad 1 \leq n \leq N-1 , 
\end{cases}
\\
E\, d_n
&=
\begin{cases}
c_{s-1}, & n = 0, \\
[n][\lambda-2s-n]\,d_{n-1}
+ c_{n+s-1},
& 1 \leq n \leq N-s, \\
[n][\lambda-2s-n]\,d_{n-1}, 
& p-s+1 \leq n \leq N-1,
\end{cases}
\\
F\, c_n
&=
\begin{cases}
c_{n+1}, & 0 \leq n \leq N-2, \\
0, & n = N-1,
\end{cases}
\qquad
F\, d_n
=
\begin{cases}
d_{n+1}, & 0 \leq n \leq N-2, \\
0, & n = N-1.
\end{cases}
\end{aligned}
$$
Then, for generic $\lambda$, ${\mathcal Y}(\lambda, s)$ is isomorphic to the direct sum ${\mathcal X}(\lambda)\oplus{\mathcal X}(\lambda-2s)$, 
where ${\mathcal X}(\lambda)$ is identified with the subspace spanned by
$\{c_0, \cdots, c_{N-1}\}$ and 
${\mathcal X}(\lambda-2s)$ is identified with the subspace spanned by
$\left\{d_{0} - \dfrac{c_{s}}{[s][\lambda-s]}, \cdots, d_{N-s-1} - \dfrac{c_{N-1}}{[s][\lambda-s]}, d_{N-s} , \cdots, d_{N-1}\right\}$.  
%
%
%
Let $\rho^\lambda(u)$ be the representation matrix of $u \in \widetilde{\mathcal U}_q(sl_2)$ on ${\mathcal X}(\lambda)$ with respect to the basis $\{v_0^\lambda, \cdots, v_{N-1}^\lambda\}$ as before and  
$\rho^{(\lambda, s)}(u)$ be the representation matrix of $u$ on ${\mathcal Y}(\lambda, s)$ with respect to the above basis
$\{c_n^{(\lambda, s)}, d_n^{(\lambda, s)}; 0 \leq n \leq N-1\}$.  
 $\rho^{\lambda}(u)_{n,n}$ and  $\rho^{\lambda-2s}(u)_{n,n}$ respectively. 
 %
 %
 %
\par
Now we express the symmetric linear function $G_s$ of the small quantum group $\overline{\mathcal U}_q(sl_2)$ by using the derivatives of the diagonal elements of the highest weight representations of the medium quantum group $\sresq$. 
\bigskip\par\noindent
{\bf Lemma 3.}
\begin{it}
$G_s(u)$ is given by
\begin{multline*}
G_s(u) = 
\dfrac{N\,\{1\}}{2 \, \pi \, i \, [s]} \, 
\\
\dfrac{d}{d\lambda}  \left.
\left(-\sum_{n=0}^{N-1}\left(\rho^{2N+\lambda-2s}(u)_{n,n}-
\rho^\lambda(u)_{n,n}\right)
+\sum_{n=0}^{N-s-1} \left(\rho^{2N+\lambda-2s}(u)_{n,n} - 
\rho^{\lambda-2s}(u)_{n,n}\right)
\right)\right|_{\lambda=s}
. 
\end{multline*}
\end{it}
\bigskip\par\noindent
{\sc Proof.}
There is a natural projection from $\widetilde{\mathcal U}_q(sl_2)$ to $\overline{\mathcal U}_q(sl_2)$ sending $K^{2N}$ to $1$,  and let $\hat u$ be a pre-image of $u$ in $\widetilde{\mathcal U}_q(sl_2)$.  
Then the action of $\rho^{(\lambda, s)}(\hat u)$ to ${\mathcal Y}(\lambda, s)$ is expressed as follows.  
$$
\begin{aligned}
\rho^{(\lambda, s)}(u)\, c_n
&=
\sum_{m=0}^{N-1}\rho^\lambda(u)_{n,m} \, c_m,
\\
\rho^{(\lambda, s)}(u)\, d_n
&=
\sum_{m=0}^{N-1}\rho^{\lambda-2s}(u)_{n,m} \, d_m + 
\sum_{m=0}^{N-1}x_{n,m}^{(\lambda, s)}(u) \, c_{m}.
\end{aligned}
$$
On the other hand, 
$$
\begin{aligned}
\rho^{(\lambda, s)}(u)
\left(d_{n} - \dfrac{c_{n+s}}{[s][\lambda-s]}\right)
&=\!
\sum_{m=0}^{N-s-1}\!
\rho^{\lambda-2s}(u)_{n,m} 
\left(d_{m} - \dfrac{c_{m+s}}{[s][\lambda-s]}\right)
+\!\!\!
\sum_{m=N-s}^{N-1}\rho^{\lambda-2s}(u)_{n,m} \,
d_{m},   
\\[-4pt]
&\qquad\qquad\qquad\qquad\qquad\qquad\qquad\qquad
 0 \leq n \leq N-s-1,
\\
\rho^{(\lambda, s)}(u)\, 
d_{n}
&=\!
\sum_{m=0}^{N-s-1}\!
\rho^{\lambda-2s}(u)_{n,m} 
\left(d_{m} - \dfrac{c_{m+s}}{[s][\lambda-s]}\right)
+\!\!\!
\sum_{m=N-s}^{N-1}\rho^{\lambda-2s}(u)_{n,m} \,
d_{m},   
\\[-4pt]
&\qquad\qquad\qquad\qquad\qquad\qquad\qquad\qquad
 N-s \leq n \leq N-1,
\end{aligned}
$$
because $d_{0} - \dfrac{c_{s}}{[s][\lambda-s]}$, $\cdots$, $d_{N-s-1} - \dfrac{c_{N-1}}{[s][\lambda-s]}$, $d_{s}$, $\cdots$, $d_{N-1}$ are identified with the weight vectors $v_0^{\lambda-s}$, $\cdots$, $v_{N-1}^{\lambda-s}$ of ${\mathcal X}(\lambda-2s)$ respectively.  
Therefore,
\begin{equation}
x_{n,m}^{(\lambda, s)}(u) = 
\begin{cases}
\dfrac{\rho^\lambda(u)_{n+s,m} - \rho^{\lambda-2s}(u)_{n,m-s}}{[s][\lambda-s]}, 
& 0 \leq n \leq N-s-1 \\
-\dfrac{\rho^{\lambda-2s}(u)_{n,m-s}}{[s][\lambda-s]}, & N-s \leq n \leq N-1, 
\end{cases}
\label{eq:x}
\end{equation}
where $\rho^{\lambda-2s}(u)_{n,m-s}$ is considered to be $0$ if $m-s < 0$.   
From \eqref{eq:x},  we get
$$
\begin{aligned}
\lim_{\lambda\to s}x_{n, n+s}^{(\lambda, s)}(u) 
&= 
\dfrac{N\, \{1\}}{2 \, \pi \, i \, [s]}  
\left.\dfrac{d}{d\lambda} 
\left(\rho^\lambda(u)_{n+s,n+s} - \rho^{\lambda-2s}(u)_{n,n}
\right)\right|_{\lambda=s},
\\
\lim_{\lambda\to 2N-s}x_{n, n+s}^{(\lambda, N-s)}(u) 
&= 
-\dfrac{N\, \{1\}}{2 \, \pi \, i \,  [s]} 
\left.
\dfrac{d}{d\lambda} \, 
\left(\rho^{\lambda}(u)_{n+N-s,n+N-s} - \rho^{\lambda-2N+2s}(u)_{n,n}
\right)\right|_{\lambda=2N-s}  
\end{aligned}
$$
for $0 \leq n \leq N-s-1$.  
Since the symmetric linear function $G_s(u)$ is given by
$$
G_s(u) = 
\sum_{n=0}^{s-1} x_{n, n+N-s}^{(2N-s, N-s)}(u)
+
\sum_{n=0}^{N-s-1} x_{n, n+s}^{(s, s)}(u),
$$
we have
$$
\begin{aligned}
G_s(u) &= 
\dfrac{N\,\{1\}}{2 \, \pi \, i \, [s]} \, \times
\\[-4pt]
&\dfrac{d}{d\lambda}  \left.
\left(-\sum_{n=0}^{s-1} \rho^{2N+\lambda-2s}(u)_{n+N-s,n+N-s}+
\sum_{n=0}^{N-1} \rho^\lambda(u)_{n,n} - 
\sum_{n=0}^{N-s-1} \rho^{\lambda-2s}(u)_{n,n}
\right)\right|_{\lambda=s}
\\
&= 
\dfrac{N\, \{1\}}{2 \, \pi \, i \, [s]} \, \times
\\[-4pt]
\dfrac{d}{d\lambda}  &\left.
\left(-\sum_{n=0}^{N-1}\left( \rho^{2N+\lambda-2s}(u)_{n,n}-
 \rho^\lambda(u)_{n,n}\right)
+\sum_{n=0}^{N-s-1} \left(\rho^{2N+\lambda-2s}(u)_{n,n} - 
\rho^{\lambda-2s}(u)_{n,n}\right)
\right)\right|_{\lambda=s}
. 
\end{aligned}
$$
 \qed
\subsection{Non-triviality of $\psi_{\phi}(U_\pm)$}
Here we show the following.  
\medskip\par\noindent
{\bf Lemma 4.} 
\begin{it}
$\left|\psi_\phi(U_\pm)\right| =1$ and is not equal to $0$.  
\end{it}
\medskip\par\noindent
{\sc Proof.}
Let $T_\pm$ be the tangle corresponding to $U_\pm$ and
$u_\pm$ be the universal invariant of $T_\pm$.  
Then, Proposition 7 in \cite{JM} shows that the scalar corresponding to
$\rho^\lambda(u_\pm)$ is $q^{\pm\frac{(\lambda-1)(\lambda+1-2N)}{2}}$.  
Let $\tilde u_\pm$ be the universal invariant for $U_\pm$, then
$\tilde u_\pm = K^{N+1} \, u_\pm$.  
We know that $T_s(\tilde u_\pm)=0$ and $\chi^\lambda(\tilde u_\pm) = 0$, 
$\phi(\tilde u_\pm)$ is computed as follows. 
$$
\begin{aligned}
\phi(\tilde u_\pm)
&=
\dfrac{1}{N\sqrt{2N}}\, 
\sum_{s=1}^{N-1} (-1)^{s-1} \, [s]^2 \, G_s(\tilde u_\pm)
\\
&=
\dfrac{1}{N\sqrt{2N}} 
\sum_{s=1}^{N-1} 
\dfrac{N\, \{s\}\, (-1)^{s-1}}{2\, \pi \, i} \, 
\dfrac{d}{d\lambda} 
\left.\left(\sum_{n=0}^{N-s-1} \left(\rho^{2N+\lambda-2s}(\tilde u_\pm)_{n,n} - 
\rho^{\lambda-2s}(\tilde u_\pm)_{n,n}\right)
\right)\right|_{\lambda=s}, 
\end{aligned}
$$
since $\chi^\lambda(K^{N+1} \, u_\pm) = 0$.  
We know that 
$$
\rho^\lambda(K^{N+1}\, u_\pm)_{nn} =  
q^{\pm\frac{(\lambda-1)(\lambda+1-2N)}{2}} \,
q^{(\lambda-2n)(N+1)}
= 
-q^{\pm\frac{\lambda^2 - 2N\lambda -1}{2}} \,
q^{\lambda-2n + N\lambda}.
$$
Hence, by using Lemma 3,  we have
\begin{multline*}
\phi(\tilde u_+)
=
\dfrac{1}{N\sqrt{2N}}\, 
\sum_{s=1}^{N-1} 
\dfrac{N \, \{1\}\, (-1)^s \, [s]}{2\, i} \, 
\sum_{n=0}^{N-s-1}q^{\frac{s^2 +2Ns -1}{2}} \,
q^{-s-2n -Ns}
\\
=
-\dfrac{\{1\} }{ 2 \,i\, \sqrt{2N}}\, 
\sum_{s=1}^{N-1} 
\, (-1)^s\, [s]^2 \, 
q^{\frac{s^2+1}{2}} 
=
\dfrac{i}{2\sqrt{2N}} \, 
g_{4N} \, q^{\frac{-N^2-N-1}{2}}
,
\end{multline*}
where $g_{4N} = \sum_{j=0}^{4N-1} q^{j^2/2} = 2\, (1+i)\, \sqrt{N}$ by Corollary 1.2.3 of \cite{BEW} and $\left|\phi(\tilde u_+)\right|=1$.  
Since $\phi(\tilde u_-)$ is the complex conjugate of $\phi(\tilde u_+)$, its absolute value is also  equal to $1$.  
\qed
\subsection{Quantum $SO(3)$ version}
For the Witten-Reshetikhin-Turaev invariant of three manifolds, its $SO(3)$ version is introduced in \cite{KM}.  
Here, we  introduce the $SO(3)$ version of the generalized Kashaev invariant.  
For the Hennings invariant,  $SO(3)$ version is already introduced in \cite{ChKuSr}.  
\par
Let $N$ be a positive odd integer, $\overline{\mathcal U}_q(so_3)$ be the quotient of $\overline{\mathcal U}_q(sl_2)$ obtained by adding a relation $K^N = 1$ to $\overline{\mathcal U}_q(sl_2)$, and $pr$ be the canonical projection from $\overline{\mathcal U}_q(sl_2)$ to $\overline{\mathcal U}_q(so_3)$.  
We also have an inclusion $\iota$ from  $\overline{\mathcal U}_q(so_3)$ to  $\overline{\mathcal U}_q(sl_2)$ sending the generators $E$, $F$, $K$ of $\overline{\mathcal U}_q(so_3)$ to $ E(K^N+1)/2$, $F(K^N+1)/2$, $K(K^N+1)/2$, which is  a Hopf algebra homomorphism since $K^N$ is a center of $\overline{\mathcal U}_q(sl_2)$.  
Then, everything we constructed before for $\overline{\mathcal U}_q(sl_2)$ can be defined for  $\overline{\mathcal U}_q(so_3)$  through $\iota$.  
Any representation of $\overline{\mathcal U}_q(sl_2)$ in \S1.5  parametrized by even $s$ maps $K^N$ to $-1$, and it maps $K^N+1$ to $0$.  
Therefore, the symmetric linear functions parametrized by even $s$ vanishes.  
Let $\phi^{SO(3)} = {\sqrt{2}}\,\phi \circ \iota$,  then
$$
\phi^{SO(3)}(u)
=
\sqrt{2} \, \left(
\alpha_N\, T_N(u) +
\sum_{s=0}^{\frac{N-3}{2}}
\left(
\alpha_{2s+1} \, T_{2s+1}(u) + \beta_{2s+1} \, G_{2s+1}(u)
\right)\right).  
$$
\par
For later use, we compute $\phi^{SO(3)}(\tilde u_\pm)$ and check their non-triviality for the universal invariant $\tilde u_\pm$ of trivial $\pm1$ framed knot $U_\pm$. 
\medskip\par\noindent
{\bf Lemma 5. }
\begin{it}
 $\left|\phi^{SO(3)}(\tilde u_\pm)\right|=1$ and is not equal to zero.  
 \end{it}
 \smallskip\par\noindent
 {\sc Proof.}
We restrict the parameter $s$ to odd integers in the computation in \S4.4.  
\begin{multline*}
\phi^{SO(3)}(\tilde u_+)
=
\dfrac{1}{N\sqrt{N}}\,\sum_{s=0}^{\frac{N-3}{2}}\dfrac{-N\, \{2s+1\} \, }{2\, i}
\sum_{n=0}^{N-2s-2} q^{\frac{(2s+1)^2 + 2N(2s+1) - 1}{2}} \, 
q^{-2s-1-2n-N(2s+1)}
\\
=
\dfrac{ \{1\}}{2\, i\, \sqrt{N}}\,
\sum_{s=0}^{\frac{N-3}{2}} [2s+1]^2\,  q^{\frac{(2s+1)^2 + 1}{2}}
=
\dfrac{\tilde g_{4N}\, q^{-1/2}}{2\, i\, \sqrt{N}},
\end{multline*}
where $\tilde g_{4N} = \sum_{s=0}^{2N-1} q^{\frac{(2s+1)^2}{2}}$.  
Since 
$$
\sum_{s=0}^{2N-1} q^{2 s^2} = 2\, g_{N} =  2\, \sqrt{N}\ \ 
(N \equiv 1 \mod 4)\ \  \text{or}\ \  2 \, i \, \sqrt{N}\ \ (N \equiv 3 \mod 4)
$$ 
by Corollary 1.2.3 of \cite{BEW}, $\left|\tilde g_{4N}\right| = \left|g_{4N}  - 2 g_{N}\right| = 2\, \sqrt{N}$, which implies that   $\left|\phi^{SO(3)}(\tilde u_+)\right|= 1$.  
Similarly, we get $\left|\phi^{SO(3)}(\tilde u_-)\right|=1$.  
\qed
\medskip\par\noindent
{\bf Theorem 6.}
\begin{it}
Let $M$ be a three manifold obtained by the surgery along a framed link $L$ in $S^3$, 
$\widetilde K$ be a link  in $M$,  and  $\widehat K$ be the pre-image of $\widetilde K$ in $S^3$.  
Let  $T$ be a tangle obtained from $\widehat K \cup L$, and $s_+(L)$, $s_-(L)$ are numbers of the positive and negative eigenvalues of the linking matrix of $L$.  
Let 
$$
\widetilde\GK_N^{SO(3)}(T) 
= 
A_{N, \cdots, N, \phi^{SO(3)},\cdots,  \phi^{SO(3)}}(T)
$$
and
$$
\GK_N^{SO(3)}(T)
=
{\psi_{\phi^{SO(3)}}(U_+)^{-s_+(L)}\, \psi_{\phi^{SO(3)}}(U_-)^{-s_-(L)}\, 
{\widetilde\GK_N^{SO(3)}(T) }}.
$$
Then $\GK_N^{SO(3)}(\widetilde K)$ is an invariant of $\widetilde K$.  
\end{it} 
\medskip
\par
We call $\GK_N^{SO(3)}(T)$ the {\it $SO(3)$ version of the generalized Kashaev invariant}.  
\subsection{Relation to other invariants}
We express the generalized Kashaev invariants  $\GK_N(\widetilde K)$ and $\GK_N^{SO(3)}(\widetilde K)$ by using the colored Alexander invariants and the colored Jones invariants.  
$\GK_N(\widetilde K)$ and $\GK_N^{SO(3)}(\widetilde K)$  are given as follows.  
$$
\widetilde\GK_N(\widetilde K) = 
\sum_\nu\Big(
\sum_{t_1, t_2, \cdots, t_p=0}^N\prod_{i=1}^p \left(\alpha_{t_i}\,T_{t_i}(u_{r+i}^\nu)+
\beta_{t_i}\,G_{t_i}(u_{r+i}^\nu)\right)\Big)  \, \prod_{j=2}^r T_N(u_j^\nu) \, u_1^\nu\, \be_N,    
$$
$$
\widetilde\GK_N^{SO(3)}(\widetilde K) = 
\sqrt{2} \, \sum_\nu
\Big(
\sum_{\begin{matrix}\scriptstyle t_1, t_2, \cdots, t_p=0\\[-2pt] \scriptstyle t_i : \text{odd}\end{matrix}}^N\prod_{i=1}^p \left(\alpha_{t_i}\,T_{t_i})(u_{r+i}^\nu)+
\beta_{t_i}\,G_{t_i}(u_{r+i}^\nu)\right)\Big)  \, \prod_{j=2}^r T_N(u_j^\nu) \, u_1^\nu\, \be_N.   
$$
\par
From now on, we consider the case that $p=1$, i.e. the framed link $L$ defining the three manifold $M$ is a knot.  
Then $\widetilde\GK_N(\widetilde K)$ and $\widetilde\GK_N^{SO(3)}(\widetilde K)$ are given by the colored Alexander invariants and the colored Jones invariants as in the following theorem.  
Here $V_{N, \cdots, N, t}(\widehat K\cup L)$ is the colored Jones invariant at $q = \exp(\pi\, i/N)$, which is normalized to be $1$ for the trivial knot.  
\smallskip\par\noindent
{\bf Theorem 7.}
\begin{it}
Let $\widetilde K$,  $M$, $L$, $\widehat K$ and $T$ be  as in Theorem 6.  
Then the generalized Kashaev invariants $\widetilde\GK_N(\widetilde K)$ and $\widetilde\GK_N^{SO(3)}(\widetilde K)$ are expressed in terms of the colored Alexander invariants, their derivatives and  the colored Jones invariants of $K \cup L$ as follows.  
\begin{multline}
\widetilde\GK_N(\widetilde K)
=
-\dfrac{1}{\sqrt{2N}} \, 
\left(\vphantom{\sum_{t=1}^{N-1}(-1)^{t}}
(-i)^{N-1}
\ADO_{N,N, \cdots, N, 0}(\widehat K \cup L)  
 +  
i^{N-1} \, \ADO_{N,N, \cdots, N, N}(\widehat K \cup L)
 \right.
\\[-4pt]
\left.
+\sum_{t=1}^{N-1}(-1)^{t} \, \{t\}_+\,{(-i)^{N-1}} \, 
\ADO_{N,N, \cdots, N, t}(\widehat K \cup L) 
+
\right.
\\
\sum_{t=1}^{N-1}(-i)^{N-1} (-1)^t  \{t\}
\left( \dfrac{N}{2 \pi i} \, 
\dfrac{d}{d\mu}\left(-\ADO_{N, \cdots, N, 2N+\mu-2t}^N(\widehat K\cup L)+\ADO_{N, \cdots, N, \mu}^N(\widehat K\cup L)
\right)\right|_{\mu=t}
\\[-4pt]
\left.
+
f
\, N \, \ADO_{N, \cdots, N, t}(\widehat K \cup L)- i^{N-1}\, f \,V_{N, \cdots, N, t}(\widehat K \cup L) 
\Big)\vphantom{\sum_{t=1}^{N-1}(-1)^{t}}\right) ,  
\label{eq:gk}
\end{multline}
\begin{multline}
\widetilde\GK_N^{SO(3)}(\widetilde K)
=
\\
\dfrac{(-1)^{\frac{N-1}{2}}}{\sqrt{N}} \, 
\left(\vphantom{\sum_{t=1}^{N-1}(-1)^{t}}
\sum_{t=0}^{\frac{N-3}{2}} \{2t+1\}_+\, \, 
\ADO_{N,N, \cdots, N, 2t+1}(\widehat K \cup L) 
-\ADO_{N,N, \cdots, N, N}(\widehat K \cup L)+
\right.
\\[-4pt]
\sum_{t=0}^{\frac{N-3}{2}} \{2t+1\}\,
\left( \dfrac{N}{2 \pi i} \, 
\dfrac{d}{d\mu}\left(-\ADO_{N, \cdots, N, 2N+\mu-4t-2}^N(\widehat K\cup L)+
\left.
\ADO_{N, \cdots, N, \mu}^N(\widehat K\cup L)
\right)\right|_{\mu=2t+1}
+
\right.
\\[-4pt]
\left.\left.
f
\, N \, \ADO_{N, \cdots, N, 2t+1}(\widehat K \cup L)-(-1)^{\frac{N-1}{2}}\, f \,V_{N, \cdots, N, 2t+1}(\widehat K \cup L) 
\vphantom{\dfrac{N}{2 \pi i}}
\right)\vphantom{\sum_{t=1}^{\frac{N-1}{2}}}\right).   
\label{eq:gkso3}
\end{multline}
\end{it}
\medskip\par\noindent
{\sc Proof. } 
Since $L$ is a knot, $\widetilde\GK_N(\widetilde K)$ and $\widetilde\GK_N^{SO(3)}(\widetilde K)$ is given as follows.  
\begin{multline*}
\widetilde\GK_N(\widetilde K)
=
-\dfrac{1}{N\sqrt{2N}} \, \left(\vphantom{\sum_{t=1}^{N-1}(-1)^{t}}a_{N, T_N, \cdots, T_N, T_0}(T) + (-1)^{N-1} \, a_{N, T_N, \cdots, T_N,  T_N}(T) + \right.
\\
\left.
\sum_{t=1}^{N-1}(-1)^{t} \, \Big(\{t\}_+\,a_{N, T_N, \cdots, T_N,  T_t}(T)+
[t]^2\,a_{N, T_N, \cdots, T_N,  G_t}(T)\Big)\right) ,  
\end{multline*}
\begin{multline*}
\widetilde\GK_N^{SO(3)}(\widetilde K)
=
\dfrac{1}{N\sqrt{N}} \, \left(\vphantom{\sum_{t=1}^{N-1}(-1)^{t}} -a_{N, T_N, \cdots, T_N,  T_N}(T) + \right.
\\
\left.
\sum_{t=0}^{\frac{N-3}{2}} \Big(\{2t+1\}_+\,a_{N, T_N, \cdots, T_N,  T_{2t+1}}(T)+
[t]^2\,a_{N, T_N, \cdots, T_N,  G_{2t+1}}(T)\Big)\right) .  
\end{multline*}
We know that
\begin{multline}
a_{N, T_N, \cdots, T_N,  T_t}(T)
=
\lim_{\lambda\to N}\dfrac{i^{N-1}\, \sin \lambda\pi}{\sin(\lambda\pi/N)} \, 
\ADO_{\lambda,N, \cdots, N, t}(\widehat K \cup L) 
\\[5pt]
= 
{(-i)^{N-1}\,N} \, 
\ADO_{N,N, \cdots, N, t}(\widehat K \cup L).  
\label{eq:Tt}
\end{multline}
By using Lemma 3, we have
\begin{multline*}
a_{N, T_N, \cdots, T_N,  G_t}(T) \, \be_N
=
\\
\dfrac{N \, \{1\}}{2\, \pi \, i \, [t]} \, 
\left(\left.
\dfrac{d}{d\mu}\left(-A_{N, N, \cdots, N, 2N+\mu - 2t}(T)+A_{N, N, \cdots, N, \mu}(T)
\right)\right|_{\mu=t}\right.
\\[5pt]
\left.\left.
+ 
\sum_\nu 
\dfrac{d}{d\mu}\left(\prod_{j=2}^r T_N(u_j^\nu)  
\sum_{n=0}^{N-t-1}\left(\rho^{2N+\mu-2t}(u_{r+1}^\nu)_{n,n} - 
\rho^{\mu-2t}(u_{r+1}^\nu)_{n,n}\right)\, u_1^\nu\right)\right|_{\mu=t}\right) \, \be_N.  
\end{multline*}
From the definition of the $R$-matrix,   
$
\prod_{j=2}^r T_N(u_j^\nu)\,\rho^{\mu}(u_{r+1}^\nu)_{n,n}\, u_1^\nu
$
has period $2N$ with respect to the parameter $\mu$ except the phase factor
$q^{f (\mu-N)^2/2}$ where $f$ is the framing of $L$.   
If the framed link $\widehat K \cup L$ is given by a link diagram with blackboard framing, then the framing 
$f$ of $L$ is given by the sum of signs of the self-crossings of $L$.  
Let $h^\nu(\mu)$ be a function of period $2N$ satisfying 
$$
\prod_{j=2}^r T_N(u_j^\nu)\,\rho^{\mu}(u_{r+1}^\nu)_{n,n}\, u_1^\nu \, \be_N 
= 
q^{f (\mu-N)^2/2} \, h^\nu(\mu)\, \be_N.
$$
Then we have
\begin{multline*}
\left.
\dfrac{d}{d\mu}\left(\prod_{j=2}^r T_N(u_j^\nu)\,\left(\rho^{2N+\mu-2t}(u_{r+1}^\nu)_{n,n} - 
\rho^{\mu-2t}(u_{r+1}^\nu)_{n,n}\right)\, u_1^\nu)\right)\right|_{\mu=t}\be_N
=
\\[5pt]
\left.
\dfrac{d}{d\mu}\left(q^{f (\mu+ N-2t)^2/2} \, h^\nu(\mu-2t)\, \be_N - 
q^{f (\mu- N-2t)^2/2} \, h^\nu(\mu-2t)\right)\right|_{\mu=t}\be_N
=
\\[5pt]
2\, f\, \pi \, i\, \prod_{j=2}^r T_N(u_j^\nu)\,\rho^{2N-t}(u_{r+1}^\nu)_{n,n}\, u_1^\nu\, \be_N,   
\end{multline*}
and 
\begin{equation}
\begin{aligned}
&a_{N, T_N, \cdots, T_N,  G_t}(T) \, \be_N
=\sum_\nu \prod_{j=2}^r T_N(u_j^\nu)\, G_t(u_{r+1}^\nu) \, u_1^\nu\, \be_N
\\&= \dfrac{N \, \{1\}}{2\, \pi \, i \, [t]}  
\left((-i)^{N-1}N  
\dfrac{d}{d\mu}\left(-\ADO_{N, \cdots, N, 2N+\mu-2t}^N(\widehat K\cup L)+\ADO_{N,\cdots, N,  \mu}^N(\widehat K\cup L)
\right)\right|_{\mu=t}
\\
&\left.\qquad\qquad\qquad\qquad\qquad\qquad\qquad\qquad
+
2\, f_2 \, \pi \, i
\sum_\nu \prod_{j=2}^r T_N(u_j^\nu)\, 
T_t^-(u_{r+1}^\nu)\, u_1^\nu\right)\, \be_N
\end{aligned}
\label{eq:Gt}
\end{equation}
$$
\begin{aligned}
&=
 \dfrac{N \, \{1\}}{2\, \pi \, i \, [t]}  
\left((-i)^{N-1} N  
\dfrac{d}{d\mu}\left(-\ADO_{N, \cdots, N, 2N+\mu-2t}^N(\widehat K\cup L)+\ADO_{N, \cdots, N, \mu}^N(\widehat K\cup L)
\right)\right|_{\mu=t}
\\
&\left.\qquad\qquad\qquad\qquad\qquad\qquad\qquad\qquad
+
2\, f \, \pi \, i
\, a_{N, T_N, \cdots, T_N, T_t^-}(T)\right)\, \be_N.  
\end{aligned}
$$
Combining \eqref{eq:Tt} and \eqref{eq:Gt}, we can express $\GK_N(\widetilde K)$ by using
$ADO_{N, \cdots, N, \mu}(\widehat K\cup L)$, its derivatives and $a_{N, T_N, \cdots, T_N, T_t^-}(T)$.   
Since 
$$
a_{N, T_N, \cdots, T_N, T_t^+}(T) + a_{N,  T_N, \cdots, T_N,T_t^-}(T) = A_{N, N, \cdots, N, t}(T) = (-i)^{N-1}\, N \, \ADO_{N, \cdots, N, t}(\widehat K \cup L)
$$ 
and $a_{N, T_N, \cdots, T_N, T_t^+}(T)$ is equal to the colored Jones invariant $V_{N, \cdots, N, t}(\widehat K\cup L)$,   
we get \eqref{eq:gk} and \eqref{eq:gkso3}.  
\qed
\section{Volume conjecture for the generalized Kashaev invariant}
The invariants $\GK_N$, $\GK_N^{SO(3)}$ are generalizations of Kashaev's invariant, and we expect that the volume conjecture proposed in \cite{Ka} and \cite{MuMu} also holds for them.  
Since $\GK_N$ may vanish for some nontrivial cases,  it is better to consider $\GK_N^{SO(3)}$.    
%
%
%
\smallskip
In the rest of the paper, we compute the invariant $\GK_N^{SO(3)}$ for some examples and check  Conjecture 3 numerically.  
\subsection{Hopf link}
\begin{figure}
$$
\text{Hopf link:} \quad \widehat K \ \raisebox{-4mm}{\epsfig{file=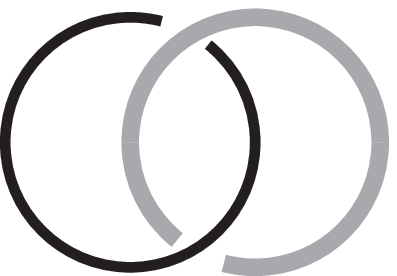, scale=0.5}}\raisebox{5mm}{$f$}  L
\qquad\qquad
\text{Whitehead link:}
\quad \raisebox{-6mm}{$\widehat K$} \ \raisebox{-10mm}{\epsfig{file=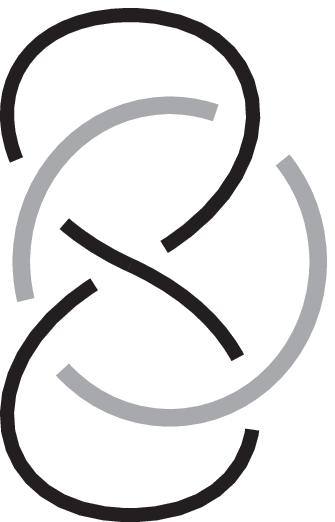, scale=0.5}}\raisebox{7mm}{$f$}  L
$$
\caption{Hopf link and Whitehead link}
\label{fig:link}
\end{figure}

Let $\widehat K\cup L$ be a Hopf link in Figure \ref{fig:link} where $f$ is the framing of $L$, and $\widetilde K$ be
the knot corresponding to $\widehat K$ in the lens space  obtained by the surgery along $L$ with framing $f$.  
We assume that the framing of $\widehat K$ is $0$.  
Let $T$ be a tangle obtained from $\widehat K\cup L$ by cutting  $\widehat K$, and $g(\mu) = (\mu^2 - 2 N \mu - 1)/2$.  
The ADO invariant of $\widehat K\cup L$ is given by 
$$
\ADO_{N, \mu}(\widehat K\cup L) =
i^{N-1} \, q^{f \, g(\mu)}, 
\qquad
\dfrac{d}{d\mu} ADO_{N, \mu}(\widehat K\cup L)
=
\dfrac{\pi \, i^N}{N}\,  (\mu-N) \, q^{f \, g(\mu)}.
$$
From the colored Jones invariant, we have 
$$
a_{N, T_t^+}(T) = 
q^{f \, g(t)} \, \lim_{\lambda \to N}\dfrac{[\lambda\, t]}{[\lambda]} =
q^{f\, g(t)} \, t  
$$
for odd $t$ ,  and  
$$
a_{N, T_t^-}(T) = 
A_{N, t}(T) - a_{N, T_t^+}(T) = 
q^{f\, g(t)} \, (N-t).  
$$
Therefore
$$
\widetilde\GK_N^{SO(3)}(\widetilde K) =
\dfrac{1}{\sqrt{N}} \, \sum_{t=0}^{N-1}
q^{2t+1} \, q^{f\, g(2t+1)}
.  
$$
This implies   $|\GK_N^{SO(3)}(\widetilde K)| \leq\sqrt{N}$ and
$$
\lim_{N\to\infty}\dfrac{2\, \pi \, \log |\GK_N^{SO(3)}(\widetilde K)|}{N}
=
0.  
$$
\subsection{Whitehead link}
We do the same thing for a Whitehead link $\widehat K \cup L$ in Figure \ref{fig:link}.  
Let  $f$ be the framing of $L$, and $\widetilde K$ be
the knot corresponding to $\widehat K$ in the lens space  obtained by the surgery along $L$ with framing $f$.  
We assume that the framing of $\widehat K$ is $0$.  
Let $T$ be a tangle of $\widehat K\cup L$ obtained by cutting a point on $\widehat K$.  
The colored Alexander invariant of $\widehat K\cup L$ is given in \cite{JM}.  
Let $\{a;, k\} = (q^a-q^{-a}) (q^{a-1}-q^{-q+1}) \cdots (q^{a-k+1} - q^{-a+k-1})$, then
$$
\ADO_{N, \mu}(\widehat K\cup L) =
\dfrac{ - i \, q^{f \, g(\mu)}}
{2 \, N \,  \sin \pi \mu} \,
\sum_{j = \left[\frac{N}{2}\right] }^{N-1}
(-1)^j \, q^{\frac{- j^2}{2}
-\frac{3 j}{2}}
\,
\dfrac{(\{j;j\})^3 \, \{\mu+j; 2j+1\}}
{\{2j+1; 2j+1-N\}} .   
$$
and, by using l'H\^opital's rule to obtain the limit $\mu\to t$ for an integer $t$, we have
$$
\ADO_{N, t}(\widehat K\cup L) =
\dfrac{ - i \, (-1)^t \, q^{f \, g(t)}}
{2 \, N \, \pi} \,
\sum_{j = \left[\frac{N}{2}\right] }^{N-1}
(-1)^j \, q^{\frac{- j^2}{2}
-\frac{3 j}{2}}
\,
\dfrac{(\{j;j\})^3 \,\left.\frac{d}{d\mu} \{\mu+j; 2j+1\}\right|_{\mu=t}}
{\{2j+1; 2j+1-N\}} .   
$$
Hence we have
$$
\begin{aligned}
a_{N, T_t}(T) &=
\dfrac{ (-i)^{N} \, (-1)^t\, q^{f \, g(t)}}{2 \, \pi} \,
\sum_{j = \left[\frac{N}{2}\right] }^{N-1}
(-1)^j \, q^{\frac{- j^2}{2}
-\frac{3 j}{2}}
\,
\dfrac{ \{j;j\}^3 \,\left.\frac{d}{d\mu} \{\mu+j; 2j+1\}\right|_{\mu=t}}
{\{2j+1; 2j+1-N\}}
\\
&=
\dfrac{ -i^{N} \, (-1)^t\, q^{f \, g(t)}}{2 \, \pi} \,
\sum_{j = \left[\frac{N}{2}\right] }^{N-1}
(-1)^j \, q^{\frac{- j^2}{2}
-\frac{3 j}{2}}
\,
\dfrac{ \{j;j\}^3 \,\left.\frac{d}{d\mu} \{\mu+j; 2j+1\}\right|_{\mu=t}}
{\{2j+1-N; 2j+1-N\}}.
\end{aligned}
$$
The colored Jones invariant $V_{N, t}(\widehat K \cup L)$  is given in \cite{H} which is redormulated  as follows for $q = e^{\pi i/N}$.  
$$
V_{N, t}(\widehat K \cup L)
=
\begin{cases}
\displaystyle
(-1)^{t-1} \, q^{f \, g(t)}\,\sum_{j=0}^{\left[\frac{N}{2}\right]-1} 
(-1)^j \, q^{-\frac{j^2}{2}-\frac{3j}{2}} \,
\frac{\{j; j\}^2 \, \{t+j; 2j+1\}}
{\{2j+1; j+1\}} \, +
\\[10pt]
\hfill
\displaystyle
\frac{ i^N \, (-1)^{t-1} \,q^{f \, g(t)}}{2 \, \pi}
\sum_{j=\left[\frac{N}{2}\right]}^{t-1} 
(-1)^j  q^{-\frac{j^2}{2}-\frac{3j}{2}} 
\frac{\{j; j\}^3 \,\left.\frac{d}{du} \{\mu+j; 2j+1\}\right|_{\mu=t}}
{ \{2j+1-N; 2j+1-N\}},
\\
\displaystyle
\hfill \text{if $t > \left[\frac{N}{2}\right]$,}
\\[20pt]
\displaystyle
(-1)^{t-1} \, q^{f \, g(t)}\,\sum_{j=0}^{t-1} 
(-1)^j \, q^{-\frac{j^2}{2}-\frac{3j}{2}} \,
\frac{\{j; j\}^2 \, \{t+j; 2j+1\}}
{\{2j+1; j+1\}},
\qquad \text{if $t \leq \left[\frac{N}{2}\right]$}. 
\end{cases}
$$
Now we check the following conjecture 
by substituting the above formulas to \eqref{eq:gkso3}.  
This conjecture is stronger than Conjecture 3. 
\medskip\par\noindent
{\bf Conjecture 4.}
\begin{it}
The $SO(3)$ version of the generalized  Kashaev invariant satisfies the following.  
$$
\lim_{N\to\infty} \pi\, \log \dfrac{ \GK_N^{SO(3)}(\widetilde K)}{\GK_{N-2}^{SO(3)}(\widetilde K)}
\equiv
\Vol(\widetilde K) + i \CS(\widetilde K)
\quad
\text{$\mod$ $\dfrac{\pi^2 i}{2} $}.   
$$  
\end{it}
\medskip\par
For the knot $\widetilde K$ in  lens spaces comes from the Whitehead link as above,   
the results  of numeric computation are exposed in Table \ref{table}.  
The framing $f$ of $L$ varies from $-5$ to $10$ and 
$N = 83$, $123$, $183$, $245$.  
For the cases $f = 0, 1, 2, 3, 4$, the knot complements are not hyperbolic.  
The volumes and Chern-Simons invariants are obtained from the software SnapPea and its cusped census table created by Jeff Weeks.  
The values seems to converge to the complex volume $\Vol(\widetilde K) + i \CS(\widetilde K) \mod \pi^2/2  i$ when $f \equiv \hspace{-4mm}  \setminus\ \ 2 \mod 4$.  
\begin{table}[thb]
\begin{tabular}{|c||c|c|c|c||c|}
\hline
$f \setminus N$ & \quad\  \ 83\quad &\quad\  123&\quad\ 183& \quad\  245 &  \small $\ \ \Vol + i \CS$\\
\hline
-5 & 
$\begin{matrix}3.52627 \\+\ 3.77047\, i\end{matrix}$ & 
$\begin{matrix}3.45119\\+\   3.77611 \, i\end{matrix}$ & 
$\begin{matrix}3.40037\\+\   3.77866\, i\end{matrix}$ & 
$\begin{matrix}3.37410\\+\   3.77958 \, i\end{matrix}$ &     
$\begin{matrix}3.29690\\-\  1.15407\, i\end{matrix}$ 
\\ \hline
-4 &$\begin{matrix}3.40671\\-\  0.97724 \, i\end{matrix}$ & 
$\begin{matrix}3.33159\\-\  0.97243\, i\end{matrix}$  & 
$\begin{matrix}3.28077\\-\  0.97025\, i\end{matrix}$ &
$\begin{matrix}3.25449\  - \\0.96946 \, i\end{matrix}$ &   
$\begin{matrix}3.17729\\-\  0.96847\, i\end{matrix}$ 
\\ \hline
-3 & 
$\begin{matrix}3.21855\\+\  4.19927\, i\end{matrix}$ & 
$\begin{matrix}3.14342 \\+\  4.20327 \, i\end{matrix}$& 
$\begin{matrix}3.09260\\+\ 4.20508\, i\end{matrix}$ &
$\begin{matrix}3.06632\\+\  4.20574 \, i\end{matrix}$  & 
$\begin{matrix}2.98912\\+\  4.20662\, i\end{matrix}$ 
\\ \hline
-2 &
$\begin{matrix}-0.20084\\-\  3.95382 \, i\end{matrix}$ & 
$\begin{matrix}0.64312\\-\  3.93690 \, i\end{matrix}$& 
$\begin{matrix}-1.15661\\-\  3.63002\, i\end{matrix}$ &
$\begin{matrix}0.30569\\-\  3.968664 \, i\end{matrix}$ &   
$\begin{matrix}2.66674\\-\ 0.41123\, i\end{matrix}$
\\ \hline
-1 & 
$\begin{matrix}2.25923\\+\  4.93040\, i\end{matrix}$ &
$\begin{matrix}2.18415\\+\  4.93281 \, i\end{matrix}$  & 
$\begin{matrix}2.13335\\+\  4.93391\, i\end{matrix}$ &
$\begin{matrix}2.10708\\+\  4.93430 \, i\end{matrix}$ &   
$\begin{matrix}2.02988\\+\  0\, i\end{matrix}$
\\ \hline
\ 0 &
$\begin{matrix}0.30651\\-\  0.00294 \, i\end{matrix}$ & 
$\begin{matrix}0.20601\\-\  0.00133 \, i\end{matrix}$&
$\begin{matrix}0.13809\\-\  0.000609 \, i\end{matrix}$ &
$\begin{matrix}0.10300\\-\  0.00033 \,i\end{matrix} $ &
\begin{tabular}{ll}non-\\hyperbolic\end{tabular}
\\ \hline
\ 1 & 
$\begin{matrix}0.22776\\+\  3.28809 \, i\end{matrix}$& 
$\begin{matrix}0.15482\\+\  3.29048 \, i\end{matrix}$ &
$\begin{matrix}0.10340\\+\  3.28898\, i\end{matrix}$ &
$\begin{matrix}0.07714\\+\  3.29014 \, i\end{matrix}$ & 
\begin{tabular}{ll}non-\\hyperbolic\end{tabular}
\\ \hline
\ 2 &
$\begin{matrix}0.23123\\-\  4.93300 \, i\end{matrix}$ &
$\begin{matrix}0.15525\\-\  4.93385 \, i\end{matrix}$ &
$\begin{matrix}0.10398\\-\  4.93430 \, i\end{matrix}$ &
$\begin{matrix}0.07752\\-\  4.93449 \, i\end{matrix}$ & 
\begin{tabular}{ll}non-\\hyperbolic\end{tabular}
\\ \hline
\ 3 &
$\begin{matrix}0.35286\\-\  1.37902 \, i\end{matrix}$ &
$\begin{matrix}0.255233\\-\  1.361564 \, i\end{matrix}$  &
$\begin{matrix}0.184937\\-\  1.344945 \, i\end{matrix}$  & 
$\begin{matrix} 0.00758\\+\  4.62530 \, i\end{matrix}$ & 
\begin{tabular}{ll}non-\\hyperbolic\end{tabular}
\\ \hline
\ 4 &
$\begin{matrix}-0.06551\\-\  4.58143 \, i\end{matrix}$ &
$\begin{matrix}-0.09561\\-\  4.64719 \, i\end{matrix}$ &
$\begin{matrix}0.31686\\+\  4.70468 \, i\end{matrix}$ &
$\begin{matrix} 0.30301\\-\  4.75063 \, i\end{matrix}$ & 
\begin{tabular}{ll}non-\\hyperbolic\end{tabular}
\\ \hline
\ 5 &
$\begin{matrix}2.25936\\+\  0.00547 \, i\end{matrix}$ &
$\begin{matrix}2.18421\\+\  0.00247 \, i\end{matrix}$ & 
$\begin{matrix}2.13337\\+\  0.00111 \, i\end{matrix}$ &
$\begin{matrix}2.10709\\+\  0.00062 \, i\end{matrix}$  &  
$\begin{matrix}2.02988\\+\  4.93480\,  i\end{matrix}$ 
\\ \hline
\ 6 & 
$\begin{matrix}-2.74804\\-\  3.36083 \, i\end{matrix}$ & 
$\begin{matrix}-7.77535\\-\  2.77035\, i\end{matrix}$ &
$\begin{matrix}-7.74245\\-\  2.25861\, i\end{matrix}$  & 
$\begin{matrix}3.02484\\+\  0.87472 \, i\end{matrix}$ &  
$\begin{matrix}2.66674\\+\  4.52357\, i\end{matrix}$
\\ \hline
\ 7 & 
$\begin{matrix}3.21829\\+\  0.73646 \, i\end{matrix}$ & 
$\begin{matrix}3.14331\\+\  0.73195 \, i\end{matrix}$  & 
$\begin{matrix}3.09255\\+\  0.72991 \,i\end{matrix}$ &
$\begin{matrix}3.06629\\+\  0.72917 \, i\end{matrix}$  & 
$\begin{matrix} 2.98912\\-\  4.20656\, i\end{matrix}$
\\ \hline
\ 8 &
$\begin{matrix}3.40638\\-\  3.95673 \, i\end{matrix}$  & 
$\begin{matrix}3.33144\\-\  3.96200 \, i\end{matrix}$ & 
$\begin{matrix}3.28071\\-\  3.96439 \,i\end{matrix}$ &
$\begin{matrix}3.25446\\-\  3.96525 \, i\end{matrix}$  &  
$\begin{matrix}3.17729\\-\  3.96634\, i\end{matrix}$
\\ \hline
\ 9 & 
$\begin{matrix}3.52592\\+\  1.16509 \, i\end{matrix}$ &
$\begin{matrix}3.45103\\+\  1.15904 \, i\end{matrix}$  &
$\begin{matrix}3.40030\\+\  1.15630 \, i\end{matrix}$  &
$\begin{matrix}3.37406\\+\  1.15531 \, i\end{matrix}$  & 
$\begin{matrix}3.29690\\-\  3.78074 \, i\end{matrix}$
\\ \hline
10 &
$\begin{matrix}2.79822\\+\  2.69441 \, i\end{matrix}$  & 
$\begin{matrix}2.72319\\+\  2.68751 \, i\end{matrix}$ &
$\begin{matrix}2.67241\\+\  2.68438 \, i\end{matrix}$  & 
$\begin{matrix}2.64615\\+\  2.68325 \, i\end{matrix} $&   
$\begin{matrix}3.37760\\+\  3.63406\, i\end{matrix}$
\\ \hline
\end{tabular}
\\[6pt]
\begin{minipage}{14cm}
\small
The last column indicate the hyperbolic volumes and the Chern-Simons invariants given by the cusped census of  SnapPea. 
\\[0pt] 
\end{minipage}
\caption{Values of $\pi\, \log \dfrac{ \GK_N^{SO(3)}(\widetilde K)}{\GK_{N-2}^{SO(3)}(\widetilde K)} \mod \pi^2 i $}
\label{table}
\end{table}

\end{document}